\documentclass[12pt]{amsart}
\usepackage{color}
\usepackage[all]{xy}
\usepackage{amssymb}
\usepackage{tikz}
\usepackage{comment}
\usepackage[mathscr]{euscript}
\usepackage{cite}
\usepackage{enumerate}
\usetikzlibrary{arrows}

\calclayout

\date{April 9, 2018}
\newtheorem{dummy}{anything}[section] 
\newtheorem{theorem}[dummy]{Theorem}
\newtheorem*{thma}{Theorem A}
\newtheorem*{thmb}{Theorem B}

\newtheorem{lemma}[dummy]{Lemma} 
\newtheorem{proposition}[dummy]{Proposition} 
\newtheorem{corollary}[dummy]{Corollary}
 
\theoremstyle{definition}
\newtheorem{definition}[dummy]{Definition}
 \newtheorem{example}[dummy]{Example}
 
 \newtheorem{remark}[dummy]{Remark}

 \newtheorem*{question}{Question}
 \newtheorem*{acknowledgement}{Acknowledgement}

\newcommand
{\eqncount}{\setcounter{equation}{\value{dummy}}%
\addtocounter{dummy}{1}}

\newcommand{\cS}{\mathcal S}

\newcommand{\cE}{\mathcal E}

\newcommand{\bC}{\mathbf C}
\newcommand{\bD}{\mathbf D}

\newcommand{\bZ}{\mathbb Z}

\newcommand{\CP}{\mathbb C P}

\newcommand{\fQ}{Q}

\newcommand{\bbL}{\mathbb L}

\newcommand{\bbR}{\mathbb R}

\newcommand{\bbZ}{\mathbb Z}

\newcommand{\sF}{\mathscr F}
\newcommand{\cy}[1]{\bbZ/{#1}}
\newcommand{\wX}{\widetilde X}
\newcommand{\wM}{\widetilde M}
\newcommand{\wP}{\widetilde P}

\newcommand{\wH}{\widetilde H}
\newcommand{\bd}{\partial}

\newcommand{\La}{\Lambda}
\newcommand{\Zpi}{\bbZ [\pi]}

\newcommand{\mmatrix}[4]{\left (\vcenter
{\xymatrix@C-2pc@R-2pc{#1&#2\\#3&#4} }
\right )}

\newcommand{\la}{\langle}
\newcommand{\ra}{\rangle}

\DeclareMathOperator{\Hom}{Hom}
\DeclareMathOperator{\wh}{Wh}

\DeclareMathOperator{\Sharp}{\#}

\DeclareMathOperator{\Res}{Res}

\DeclareMathOperator{\Ext}{Ext}

\DeclareMathOperator{\sign}{sign}

\DeclareMathOperator{\cok}{coker}

\DeclareMathOperator{\rk}{rk}

\DeclareMathOperator{\ad}{ad}
\DeclareMathOperator{\pt}{pt}
\DeclareMathOperator{\lk}{Lk}
\DeclareMathOperator{\coker}{coker}

\newcommand{\nr}[1]{\medskip\noindent{(\textbf #1).}}
\newcommand{\nrb}[1]{\medskip\noindent{(\textbf {#1}).} }
\newcommand{\Bpi}{K(\pi,1)}

\newcommand{\RAAG}{right-angled Artin group} 

\newcommand{\RAAGs}{right-angled Artin groups} 
\newcommand{\RRAAGs}{Right-angled Artin groups} 
\newcommand{\LM}{L_M}
\DeclareMathOperator{\cd}{cd}

\newcommand{\spinp}{spin$^{+}$}
\DeclareMathOperator{\tr}{tr}
\DeclareMathOperator{\id}{id}
\begin{document}

\title[$4$-manifolds with right-angled Artin fundamental groups]
{Topological $4$-manifolds with right-angled Artin fundamental groups}
\author{Ian Hambleton and Alyson Hildum} 
\address{Department of Mathematics \& Statistics
 \newline\indent
McMaster University
 \newline\indent
Hamilton, ON  L8S 4K1, Canada}
\email{hambleton{@}mcmaster.ca}
\address{Department of Mathematics \& Computer Science
 \newline\indent
Wesleyan University
 \newline\indent
Middletown, CT 06459}
\email{ahildum@wesleyan.edu}
\thanks{Research partially supported by NSERC Discovery Grant A4000.}

\begin{abstract}\noindent
We classify closed, \spinp, topological  $4$-manifolds with fundamental group $\pi$ of cohomological dimension $\leq 3$ (up to $s$-cobordism), after stabilization by connected sum with at most $b_3(\pi)$ copies of $S^2\times S^2$. In general we  must also assume that $\pi$ also satisfies certain $K$-theory and assembly map conditions. Examples for which these conditions hold include the torsion-free fundamental groups of $3$-manifolds and  all \RAAGs\ whose defining graphs have no $4$-cliques. 
\end{abstract}
\maketitle

\section{Introduction}\label{sec:one}
Freedman \cite{Freedman82} classified simply connected, topological $4$-manifolds up to homeomorphism, and established a framework for studying 
non-simply connected $4$-manifolds. 

For a non-simply connected $4$-manifold $M$, the basic homotopy invariants are the fundamental group $\pi :=\pi_1(M)$, the second homotopy group $\pi_2(M)$, the equivariant intersection form $s_M$, and the first $k$-invariant, $k_M \in H^3(\pi;\pi_2(M))$.  These invariants give the \emph{quadratic $2$-type}
$$ \fQ(M): = [\pi_1(M), \pi_2(M), k_M, s_M]$$
which has been shown to determine the classification up to $s$-cobordism of  TOP 4-manifolds with geometrically $2$-dimensional fundamental groups (see \cite{HKT09}). For manifolds with finite fundamental groups, it is likely that additional invariants are needed (see \cite{HambletonKreck88}).

\begin{question} Are closed, oriented, spin $4$-manifolds with isometric quadratic $2$-types and torsion-free  fundamental groups always $s$-cobordant~?
\end{question}

In this paper, we study spin $4$-manifolds with fundamental groups belonging to the interesting class of \emph{\RAAGs} (or RAAGs). We build on the methods of \cite{HKT09} and \cite{Kreck99}, but new difficulties appear since our fundamental groups have cohomological dimension $> 2$ in general. Recall that a \RAAG\ is defined by a presentation associated to a finite graph $\Gamma$ (see Section \ref{RAAG section}).  An $i$-\emph{clique} in $\Gamma$ is a complete subgraph of $\Gamma$ with $i$ vertices, and we let $b_i(\pi)$ denote the number of $i$-cliques in $\Gamma$.

\begin{thma} Let $\pi$ be a \RAAG\ defined by a graph $\Gamma$ with no $4$-cliques.
 Suppose that $M$ and $N$ are closed,
  \spinp, topological $4$-manifolds with fundamental group $\pi$. Then any isometry between the quadratic $2$-types of $M$ and $N$ is stably realized by an $s$-cobordism between $M \Sharp r(S^2 \times S^2)$ and $N \Sharp r(S^2 \times S^2)$, whenever $r \geq b_3(\pi)$.
\end{thma} 
A \textbf{\spinp manifold} is a spin manifold with an additional evenness property for its \emph{dual} intersection form $h_M$  on $ \pi_2(M)^*$ (see Lemma \ref{s_M star}(ii) and Definition \ref{def:reduced SW}).
We will recall below the notion of (stable) \emph{isometry} for quadratic $2$-types.

Theorem A is an application of our main result, Theorem \ref{thm:maintame},  which covers the large class of fundamental groups with cohomological dimension  $\cd \pi \leq 3$  satisfying certain assembly map conditions  (see the properties (W-AA) given in Definition \ref{def:WAA}). For example, this class contains all torsion-free fundamental groups of $3$-manifolds \cite{Bartels:2014}. 
Theorem \ref{thm:maintame}  extends the results of \cite{HKT09} which handled the class of geometrically 2-dimensional groups.
If the fundamental group $\pi$ happens to be a ``good" group for topological surgery \cite{Freedman:1990}, \cite{Freedman:1995a},  \cite{Freedman:1995}, then we show that any isometry can be realized by a homeomorphism. 

 For a \RAAG\ $\pi$ these assembly conditions (W-AA) also hold \cite{Bartels:2012}, and  $\pi$ has  cohomological dimension $\cd\pi \leq 3$ if and only if the defining graph has no $4$-cliques (or equivalently if $H_4(\pi;\bZ) = 0$).  
 Thus Theorem A applies to 
 an infinite number of \RAAGs, 
 including $\pi = \bbZ^3$ (see the examples following Proposition \ref{prop:tame examples}).

 \smallskip
Here is a brief summary of the definitions in \cite{HambletonKreck88} and \cite{HKT09}. 

\begin{definition}
For an oriented 4-manifold $M$, the \emph{equivariant intersection form} is the triple $(\pi_1(M,x_0), \pi_2(M,x_0), s_M)$, where $x_0\in M$ is a base point and 
\begin{equation*}
s_M\colon  \pi_2(M,x_0) \otimes_\bbZ \pi_2(M,x_0) \to \La,
\end{equation*}
where $\La :=\bbZ[\pi_1(M,x_0)]$. This pairing is derived from the cup product on $H^2_c(\wM;\bbZ)$, where $\wM$ is the universal cover of $M$;  we identify $H^2_c(\wM;\bbZ)$ with $\pi_2(M)$ via Poincar\'e duality and the Hurewicz Theorem, and so $s_M$ is defined by
\begin{equation*}
s_M(x,y) = \sum_{g \in \pi} \varepsilon_0(\tilde x \cup \tilde yg^{-1}) \cdot g \in \Zpi ,
\end{equation*}
where $\tilde x, \tilde y \in H^2_c(\wM;\bbZ)$ are the images of $x,y \in \pi_2(M)$ under the composite isomorphism $\pi_2(M) \to H_2(\wM;\bbZ) \to H^2_c(\wM;\bbZ)$ and $\varepsilon_0$ is given by $\varepsilon_0\colon H^4_c(\wM;\bbZ) \to H_0(\wM;\bbZ)=\bbZ$.
\end{definition}
Unless otherwise mentioned, our  modules are  \emph{right} $\La$-modules.
This pairing is $\Lambda$-hermitian, in the sense that for all $\lambda\in\Lambda$, we have
\[ s_M( x, y\cdot \lambda)= s_M(x,y)\cdot \lambda  \quad \text{ and } \quad s_M(y, x)= \overline{s_M(x,y)},
\]
where $\lambda\mapsto \bar\lambda$ is the involution on $\Lambda$ given by the orientation character of $M$. This involution is determined by $\bar g = g^{-1}$ for $g \in \pi_1(M,x_0)$. For later reference, we note that when $M$ is spin the  term
$\varepsilon_0(\tilde x, \tilde y) \equiv 0 \pmod 2$, so $s_M$ is an \emph{even} hermitian form.

\begin{definition}\label{isometry}
An \emph{isometry} between quadratic $2$-types $\fQ(M)$ and $\fQ(M')$ is a pair $(\alpha, \beta)$, where $\alpha\colon \pi_1(M,x_0) \to \pi_1(M',x'_0)$ is an isomorphism of fundamental groups and  $$\beta\colon (\pi_2(M,x_0), s_M) \to (\pi_2(M',x'_0), s_{M'})$$ is an $\alpha$-invariant isometry of the equivariant intersection forms, such that $(\alpha^*, \beta_*^{-1})(k_{M'}) = k_{M}$.   In addition, the following diagram
$$
\vcenter{\xymatrix{0 \ar[r]& H^2(\pi; \La) \ar[r]\ar@{=}[d]& H^2(M'; \La) \ar[r]^(0.4){e_{M'}}\ar[d]_{\cong}^{\beta^*}&
\Hom_\La(H_2(M';\La),\La) \ar[r]\ar[d]_{\cong}^{\beta^*}&H^3(\pi;\La) \ar[r]\ar@{=}[d]& 0\cr
0 \ar[r]& H^2(\pi; \La) \ar[r]& H^2(M; \La) \ar[r]^(0.4){e_{M}}&
\Hom_\La(H_2(M;\La),\La) \ar[r]&H^3(\pi;\La) \ar[r]& 0}}
$$
arising from the universal coefficient spectral sequence \emph{commutes}, with maps $e_M$, $e_{M'}$ induced by evaluation, and  $\beta$ after identifying $\pi:= \pi_1(M,x_0) \cong \pi_1(M',x'_0)$ via $\alpha$. We will  assume throughout that our manifolds are connected, so that a change of base points leads to isometric intersection forms.
By a \emph{stable} isometry, we mean an isometry of quadratic $2$-types after adding a hyperbolic form $H(\La^r)$ to both sides. 
\end{definition}

Stabilization of quadratic $2$-types corresponds geometrically to connected sums with finitely many copies of $S^2 \times S^2$. Two closed topological (or smooth) $4$-manifolds $M$ and $N$ are said to be \emph{stably homeomorphic} (or \emph{stably diffeomorphic}) if 
 $$M \Sharp r(S^2 \times S^2) \approx N \Sharp s(S^2 \times S^2)$$
 are homeomorphic (or diffeomorphic) after stabilization, for some $r, s \geq 0$. Manifolds which are $s$-cobordant are stably homeomorphic, but it turns out that the $k$-invariant is not needed for the stable classification.
 
 Our main result implies the following (compare \cite[Lemma 4.5]{HKT09}, \cite[Theorem 1.5]{Kasprowski:2017}): 
 
 \begin{thmb} Let $\pi$ be a finitely presented group with  $\cd \pi \leq 3$ satisfying the assembly property \textup{(W-AA)(iii)}. If $M$ and $N$ are closed, oriented, \spinp, topological (respectively, smooth) $4$-manifolds with fundamental group $\pi$, then $M$ is stably homeomorphic (respectively, diffeomorphic) to $N$ if and only if  their equivariant intersection forms $s_M$ and $s_N$ are stably isometric.
\end{thmb}

\begin{remark} If $\pi$ is the fundamental group of a closed, oriented aspherical $3$-manifold, then Theorem B  was proved by Kasprowski, Land, Powell and Teichner \cite[Theorem 1.5]{Kasprowski:2017}, as a special case of their  stable classification of closed topological $4$-manifolds with these fundamental groups. In this setting, the authors also gave explicit invariants to completely describe the stable classification for every $w_2$-type. 
\end{remark}
\medskip
Here is a short outline of the paper.
In Section \ref{RAAG section} we recall some background on the cohomology of  right-angled Artin groups. In Section 
\ref{sec:tame} we discuss the notion of \emph{tame cohomology} for discrete groups, which is motivated by studying $4$-manifolds and may have some independent interest (see Definition \ref{def: tame}). Section  \ref{sec:fpgroups} establishes the key cohomological properties we need for groups with $\cd \pi \leq 3$, and Section \ref{sec:four} shows that the $\La$-dual $\pi_2(M)^* := \Hom_\La(\pi_2(M),  \La)$ is a stably free $\La$-module, for every closed spin $4$-manifold $M$ with $\cd \pi_1(M) \leq 3$. These results are applied in Section \ref{sec:five} to define the  \emph{reduced equivariant intersection form} on $\pi_2(M)^*$, and prove that it is non-singular. Examples are provided in Corollary \ref{cor:allgroups} of \spinp manifolds for every such fundamental group. 

 Our strategy for proving Theorem A uses Kreck's \emph{modified surgery} \cite{Kreck99} (see \cite[\S 2]{HKT09} for an overview of this theory).  In this approach one studies the bordism groups $\Omega_4(B)$ for a fibration $B \to BSTOP$. If two manifolds $M$ and $N$ admit suitable $B$-bordant reference maps, then $M$ and $N$ will be stably homeomorphic after connected sum with finitely many copies of $S^2 \times S^2$. Moreover, $M$ and $N$ will be $s$-cobordant provided a certain assembly map obstruction vanishes. 
 
 In Section \ref{sec:seven} we define an appropriate fibration $B:= B(M) = P \times BTOPSPIN$ for our setting, where $P$ is the reduced $2$-type of $M$ (see  Definition \ref{def:reduced 2type}). We show that  two manifolds $M$ and $N$ with isometric quadratic $2$-types admit reference maps into $B$.  In Sections \ref{sec:eight},  \ref{sec:nine} and \ref{sec:ten} we compute the bordism groups $\Omega_4(B) =\Omega^{Spin}_4(P)$ and give explicit criteria for deciding when two \spinp manifolds $M$ and $N$ are $B$-bordant. Finally, in Section \ref{sec:eleven} we put these ingredients together to prove our main result, Theorem \ref{thm:maintame}, and its application to Theorem A
 and Theorem B. 
  The ``extra"  stabilization needed in our setting to obtain an $s$-cobordism classification arises from the fact that our reference maps are not $2$-connected unless $H^3(\pi;\La) = 0$. When this cohomology group is non-zero, we must stabilize by $b_3(\pi)$ copies of $S^2\times S^2$ to complete the argument. 
  
  Some background material on homological algebra and the universal coefficient spectral sequence is provided for the reader's convenience in an Appendix (see Section \ref{sec:append}).

\begin{remark} We restrict to \emph{\spinp} $4$-manifolds for simplicity, but expect that analogous results hold by including the full $w_2$-type in the data (see \cite{HKT09}). To shorten notation in later sections, we will let  $\Omega^{Spin}_*(P)$ denote either the \emph{smooth} or the \emph{topological} spin bordism group of a space $P$, depending on the context.  Note that the Kirby-Siebenmann invariant \cite{kirby-siebenmann1} is determined by $\sign(M) \pmod{16}$ for spin manifolds.
\end{remark}

\begin{acknowledgement}
The authors would like to thank Diarmuid Crowley, Matthias Kreck, and Peter Teichner for useful conversations, and for pointing out mistakes in earlier versions of this work. We also thank the referees for their critical reading, and Erg\"un Yal\c c\i n for helpful information about the hypercohomology spectral sequence.
\end{acknowledgement}

\section{Cohomology of \RAAGs} \label{RAAG section}

A \emph{\RAAG} $\pi$ is a finitely generated group whose relators consist solely of commutators between generators. \RRAAGs\ are  also called \emph{graph groups} since each generator of $\pi$ can be represented by a vertex of a graph $\Gamma = \Gamma(\pi)$,  and pairs of commuting generators in $\pi$ are represented by edges in $\Gamma$ between the corresponding vertices. 

If $\pi$ has a presentation with $g$ generators and $r$ relators, we construct a handlebody with fundamental group $\pi$ using one 0-handle, $g$ 1-handles, and $r$ 2-handles attached to reflect the relations of $\pi$. In the case that $\pi$ is a \RAAG, this handlebody is homotopy equivalent to the 2-skeleton $K$ of a standard classifying space for $\pi$, known as the \emph{Salvetti complex} (see Charney \cite[\S 3.6]{Charney07}).

The integral homology and cohomology ring of $\pi$ are calculated in \cite{KimRoush80} and \cite{Charney:1995}. In \cite{Charney:1995}, it is shown that the $i$-th homology group and $i$-th cohomology group are both isomorphic to the $i$-th group of cellular chains of the Salvetti complex, and thus are  free abelian groups. In fact, the rank $b_i(\pi)$  of $H_i(\pi;\bZ)$ is equal to the number of $i$-cliques (complete subgraphs on $i$ vertices) in the defining graph $\Gamma$ for $\pi$. Thus the (co)homological dimension $\cd\pi$ of a \RAAG\ $\pi$ equals the maximum number of $i$-cliques in $\Gamma$.

In \cite{JensenMeier05}, Jensen and Meier calculate the cohomology of \RAAGs\ with \emph{group ring coefficients} using a simplicial complex $\hat\Gamma$ induced from the defining graph $\Gamma$ of $\pi$. In \cite{DavisOkun12}, Davis and Okun give a different formulation of the same theorem. We first give some necessary definitions before the statement of their results.

Let $\Gamma$ be a simplicial graph. Then the \emph{flag} complex $\hat\Gamma$ generated by $\Gamma$ is the minimal  simplicial complex in which every complete subgraph in $\Gamma$ spans a simplex. In other words, the $(i-1)$-simplices of $\hat\Gamma$ are the $i$-cliques of $\Gamma$. The \emph{link} $\lk(\sigma)$ of a simplex $\sigma$ in $\hat\Gamma$ is the collection of simplices  $\tau \in \hat\Gamma$ disjoint from $\sigma$, such that $\sigma$ and $\tau$ are sub-complexes of a higher dimensional simplex in $\hat\Gamma$. By definition, the link of the empty simplex is the entire flag complex $\hat\Gamma$ and by convention, $\dim \emptyset = -1$.  

\begin{definition}
For a simplex $\sigma\in \hat\Gamma$, we define the subgroup $\pi_\sigma \leq \pi$ to be the \RAAG\  generated by the subgraph of $\Gamma$ spanned by the vertices of $\sigma$.  By convention, $\pi_\emptyset = 1$.
\end{definition}
In general,  an induced subgraph of $\Gamma$ defines a subgroup of $\pi$ which is also a \RAAG,  but in this case,  since simplices in $\hat\Gamma$ are in bijection with  complete subgraphs in $\Gamma$, we see that $\pi_\sigma$ is a free abelian group of rank equal to $\dim\sigma +1$.  

\begin{theorem}[{Jensen-Meier \cite{JensenMeier05}, Davis-Okun \cite[Theorem 3.3]{DavisOkun12}}] \label{group ring coefficients}
Let $\Gamma$ be the defining graph for a \RAAG\ $\pi$, and let $\cS$ be the set of simplices in the induced flag complex $\hat\Gamma$. There is a spectral sequence converging to $H^*(\pi;\La)$ whose associated graded groups are given by 
\begin{equation*}
\operatorname{Gr }H^*(\pi;\La) = \bigoplus_{\sigma \in \cS} \left( \wH^{*-\dim\sigma-2}(\lk(\sigma)) \otimes \bbZ[\pi/\pi_\sigma]\right)
\end{equation*}
when $\hat\Gamma$ is not a single simplex. If $\hat\Gamma$ is a single simplex, then $\pi\cong \bbZ^n$, and $H^n(\pi;\La)= \bbZ$ is the only non-vanishing cohomology group.
\end{theorem}

\begin{example} \label{H 2 filtration}
We describe the associated graded group for $H^2(\pi;\La)$ using
Theorem \ref{group ring coefficients}. In the corresponding filtration, the empty simplex is the bottom of the filtration,  so we have the following filtration subgroups (indexed so that the top index matches the cohomology degree in question):
\begin{equation*}
\xymatrix{0\ \ar@{^(->}[r] & \sF_{0} \ar@{^(->}[r] \ar@{=}[d] & \sF_1 \ar@{^(->}[r] \ar@{->>}[d] & \sF_2=H^2(\pi;\La) \ar@{->>}[d] \\ & \sF_{0} & \sF_1/\sF_{0} & \sF_2/\sF_1} 
\end{equation*}
Note that $\wH^{-1}(\emptyset)=\bZ$ is the only non-trivial cohomology group of the empty simplex.
The filtration quotients are given by:
\begin{enumerate}\addtolength{\itemsep}{0.3\baselineskip}
\item $\sF_{0}\cong H^1(\hat\Gamma)\otimes \bbZ[\pi]$, 
\item $\sF_1/\sF_{0}\cong \bigoplus_{\sigma \in Vert(\hat\Gamma)} (\wH^0(\lk(\sigma)) \otimes \bbZ[\pi/\pi_\sigma])$, and 
\item $\sF_2/\sF_1 \cong \bigoplus_{\sigma \in Edge(\hat\Gamma)} (\wH^{-1}(\lk(\sigma)) \otimes \bbZ[\pi/\pi_\sigma])$. 
\end{enumerate}
\smallskip
The $\La$-action on the right-hand side is given by the identity on the integral cohomology of the links, 
tensored over $\bZ$ with the natural action on the induced modules $\bbZ[\pi/\pi_\sigma]$.
\end{example}

\section{\RAAGs\ with tame cohomology}
\label{sec:tame}
In this section we investigate   certain group cohomological conditions arising in the study of $4$-manifolds. Recall that finitely presented groups have either 1, 2, or infinitely many ends, where the number of ends of a group $G$ is equal to the number of ends of any Cayley graph of $G$.
  In \cite{Stallings:1968},  Stallings proved that a finitely presented group $G$ has more than one end if and only if it decomposes as a non-trivial amalgamated product or an HNN extension over a finite subgroup. 

Since  a \RAAG\  $\pi$  is infinite and torsion-free, $\pi$ will have more than one end if and only if it decomposes as a free product (see Dunwoody \cite{Dunwoody:1985}). Thus \RAAGs\ are 1-ended, or equivalently  $H^1(\pi;\La)=0$, if and only if $\pi\neq \bZ$ and the defining graphs are  connected  (see  \cite[\S 3.7]{Charney07}).
 Note that the group cohomology $H^i(\pi; \La)$ of any group $\pi$ with group ring coefficients is also a $\La$-module since $\La $ is a $\La$-$\La$ bimodule. 
 We say that $A$ is a \emph{torsion}  $\La$-module whenever $\Hom_\La(A, \La) = 0$.
 
\begin{definition}\label{def: tame}
 A finitely presented group $\pi$  has \emph{tame cohomology} if the following conditions hold: 
\begin{enumerate}
\item $\Hom_\La(H^2(\pi;\La),\La) = 0$ 
\item $\Hom_\La(H^3(\pi;\La),\La)= 0$ 
\item $\Ext^1_\La(H^3(\pi;\La),\La)=0$
\end{enumerate}
\end{definition}

 In Section \ref{sec:five} we show that when the fundamental group $\pi_1(M)$ of a closed $4$-manifold $M$ has tame cohomology,  then the dual of its equivariant intersection form $s_M$ is non-singular. We expect this property to have a key role in extending our classification results to \RAAGs\ of higher cohomological dimension.
 
  In the special case when $\cd \pi_1(M) \leq 3$, our classification results do not assume tame cohomology: we will establish conditions (ii) and (iii) in Proposition \ref{almost tame}, and prove the non-singularity of this form in Lemma \ref{s_M star}, without assuming condition (i).
 
\begin{remark}\label{tame restrictions} In the first version of this paper, we asked whether all \RAAGs\ have tame cohomology. In response,  
 Jonathan Hillman kindly provided the following example:  $H^2(\pi;\La)\cong \La$, for $\pi = (\bbZ^2*\bbZ^2)^2$.   However, we do not  know the answer for  right-angled Artin groups with no 4-cliques in the associated graph. 
 \end{remark} 

We now apply the results of Jensen-Meier and Davis-Okun cited in Theorem \ref{group ring coefficients}
 to study the  low-dimensional cohomology of \RAAGs.
   For any $\La$-module  $A$, 
we use the notation   $A^*=\Hom_\La(A,\La)$ for the \emph{dual module}.
We begin with the following observation. Let $\pi_{n,m}$ denote the free product of $n$ copies of $\bZ^2$ and $m$ copies of $\bZ$.
\begin{lemma}\label{lem:disjoint}
If $\pi$ is a \RAAG, then there exists a subgroup $\pi_{n.m} \leq \pi$ such that $H^1(\pi;\La)\cong H^1(\pi_{n.m};\Res\La)$ under restriction. If $\pi$ is a finitely generated free group, then $H^1(\pi;\La)^* = 0$.
\end{lemma}
\begin{proof}
If the defining graph $\Gamma =\Gamma(\pi)$ is connected, then 
$H^1(\pi;\La) =0$ and there is nothing to prove (take the empty free product).
Otherwise, we can write $\Gamma $ as a disjoint union of non-empty connected graphs, and let $\Gamma(n,m)\subset \Gamma$ denote a subgraph consisting of all the $m$ singleton vertex components together with one edge from each of the other $n$ connected components. The \RAAG\  defined by $\Gamma(n.m)$ is isomorphic to $\pi_{n.m}$.
Since the filtration terms for computing $H^1(\pi;\La)$ consist of $\sF_0 \cong \widetilde H^0(\hat\Gamma) \otimes \Zpi$ and 
$$\sF_1/\sF_0 \cong  \bigoplus_{\sigma \in Vert(\hat\Gamma)} (\wH^{-1}(\lk(\sigma)) \otimes \bbZ[\pi/\pi_\sigma]),
$$
it follows that $H^1(\pi;\La)$ is mapped isomorphically under restriction to $H^1(\pi_{n.m};\Res\La)$.

Suppose now that $n=0$ so that $\pi$ is a free group of rank $m$. Let $N$ be the connected sum of $m$ copies of $S^1 \times S^3$. An easy argument using the universal coefficient spectral sequence now shows that $H^1(\pi;\La)^* = 0$.
\end{proof}

The next result gives an inductive criterion for the cohomology groups to be torsion $\La$-modules.
\begin{lemma} \label{Hom duals are 0}
If $H^i(\hat\Gamma)=0$  then $\Hom_\La(H^{i+1}(\pi;\La), \La)=0$.
\end{lemma}

\begin{proof} The case $i=0$ is trivial, so 
consider the case when $i=1$. The filtrations for $H^2(\pi;\La)$ are given in Example \ref{H 2 filtration}, with $\sF_0=0$ by the assumption that $H^1(\hat\Gamma)=0$.  Consequently, $\sF_1 = \sF_1/\sF_0$. By dualizing  the short exact sequence 
$$0\to \sF_1/\sF_0 \to \sF_2 \to \sF_2/\sF_1 \to 0,$$
we get the exact sequence 
$$0 \to (\sF_2/\sF_1)^* \to (\sF_2)^* \to (\sF_1/\sF_0)^* \to \Ext_{\La}^1(\sF_2/\sF_1,\La) \to \dotsm.$$
We claim $(\sF_1/\sF_0)^*=0$ and $(\sF_2/\sF_1)^*=0$: Hom splits over a finite direct sum, and any non-zero summands in $\sF_1/\sF_0$ or $\sF_2/\sF_1$ will have torsion elements in the tensor product, and are killed by the Hom functor. For example,

\begin{eqnarray*}
(\sF_1/\sF_0)^* &=& \Hom_\La\left(\bigoplus\nolimits_{\sigma \in Vert(\hat\Gamma)} \wH^0(\lk(\sigma)) \otimes \bbZ[\pi/\pi_\sigma], \La \right) \\
&=& \bigoplus\nolimits_{\sigma \in Vert(\hat\Gamma)} \Hom_\La(\wH^0(\lk(\sigma)) \otimes \bbZ[\pi/\pi_\sigma], \La),
\end{eqnarray*}
where the direct sum is taken over all vertices $\sigma$ in $\hat\Gamma$.
By  \cite[Corollary 2.8.4]{benson1} (the Eckmann-Shapiro Lemma),
\eqncount
\begin{equation} \label{shapiro}
\Hom_{\Zpi}(\bbZ[\pi/\pi_\sigma], \Zpi) \cong \Hom_{\bbZ[\pi_\sigma]}(\bbZ, \Res_{\pi_\sigma}^{\pi}(\Zpi)) 
\end{equation}
is zero as long as $\pi_\sigma\neq 1$, or equivalently, $\sigma$ is not the empty simplex.
Thus from the long exact sequence, we see $(\sF_2)^* = H^2(\pi;\La)^* = 0$.

For any $i$, the assumption that $H^i(\hat\Gamma)=0$ implies that the bottom filtration term $\sF_0=0$. By the argument used in (\ref{shapiro}), the subquotients $\sF_k/\sF_{k-1}$, $1\leq k \leq i+1$,   in the filtration for $H^{i+1}(\pi;\La)$  are all torsion modules.
\begin{equation*}
\xymatrix{0=\sF_0 \ar@{^(->}[r] & \sF_1 \ar@{^(->}[r] \ar@{->>}[d] & \ \dotsm \ \ar@{^(->}[r] & \sF_{i-1} \ar@{^(->}[r] \ar@{->>}[d] & \sF_i \ar@{^(->}[r] \ar@{->>}[d] & \sF_{i+1}=H^{i+1}(\pi;\La) \ar@{->>}[d] \\  & \sF_1/\sF_0 & & \sF_{i-1}/\sF_{i-2} & \sF_i/\sF_{i-1} & \sF_{i+1}/\sF_i} 
\end{equation*}
By dualizing the short exact sequences
$$ 0 \to \sF_{k-1} \to \sF_k \to\sF_k/\sF_{k-1}\to 0,$$
   we see that $(\sF_{i+1})^* =H^{i+1}(\pi;\La)^* = 0$.
  \end{proof}

\begin{lemma} \label{tame}
Let $\pi$ be a \RAAG\ with defining graph $\Gamma$.
 Suppose that the induced flag complex $\hat\Gamma$  is 2-connected and that $H^1(\lk(\sigma))$ is zero for every vertex $\sigma$. Then $\pi$ has tame cohomology. 
\end{lemma} 

\begin{proof}
Since $\hat\Gamma$ is 2-connected, $H^i(\hat\Gamma)=0$ for $i=1,2$. By Lemma \ref{Hom duals are 0}, $H^{i+1}(\pi;\La)^*=0$ for $i=2,3$. Thus the first two conditions for $\pi$ to have tame cohomology are satisfied. For the last condition, we claim that $\Ext_{\La}^1(H^3(\pi;\La),\La)=0$. The filtration quotients are described below (with $\sF_3 = H^3(\pi;\La))$:
\begin{enumerate}
\item $\sF_0=H^2(\hat\Gamma)\otimes \bbZ[\pi]$,
\item $\sF_1/\sF_0$ is a direct sum of $H^1(\lk(\sigma)) \otimes \bbZ[\pi/\pi_\sigma]$ over vertices,
\item $\sF_2/\sF_1$ is a direct sum of $\wH^{0}(\lk(\sigma)) \otimes \bbZ[\pi/\pi_\sigma]$ over edges, and
\item $\sF_3/\sF_2$ is a direct sum of $\wH^{-1}(\lk(\sigma)) \otimes \bbZ[\pi/\pi_\sigma]$ over faces.
\end{enumerate}

Consider the dual of the the short exact sequence $0\to \sF_2 \to \sF_3 \to \sF_3/\sF_2 \to 0$, recalling that $(\sF_2)^*=0$:
\begin{equation*}
0 \to \Ext_{\La}^1(\sF_3/\sF_2,\La) \to \Ext_{\La}^1(\sF_3,\La) \to \Ext_{\La}^1(\sF_2,\La) \to \dotsm
\end{equation*}
Since $H^1(\lk(\sigma))=0$ for every vertex $\sigma$, $\sF_1/\sF_0=0$.  This implies that $\sF_1=0$ and so $\sF_2 \cong \sF_2/\sF_1$. The functor $\Ext_{\La}^1$ splits over a finite direct sum, so the first and third non-zero terms in the above sequence are of the form

\begin{equation*}
\Ext_{\La}^1(\sF_{i+1}/\sF_i,\La) = \bigoplus \wH^{1-i}(\lk(\sigma))
 \otimes \Ext_{\La}^1(\bbZ[\pi/\pi_\sigma],\La),
\end{equation*}
and we claim that the Ext term on the right-hand side is zero, hence $\Ext_{\La}^1(\sF_3,\La)=0$. The simplices in question are either faces (involved in the first term of the sequence) or edges (involved in the third term). By the Eckmann-Shapiro Lemma, 
$$\Ext^1_{\bbZ[\pi]}(\bbZ[\pi/\pi_\sigma],\bbZ[\pi]) \cong \Ext^1_{\bbZ[\pi_\sigma]}(\bbZ, \operatorname{Res}_{\pi_\sigma}^\pi \bbZ[\pi]).$$
 The restriction of $\bbZ[\pi]$ to $\pi_\sigma$ is just an infinite direct sum of $\bbZ[\pi_\sigma]$'s,  and
\begin{equation*}
\Ext^1_{\bbZ[\pi_\sigma]}(\bbZ, \bbZ[\pi_\sigma]) \cong H^1(\pi_\sigma; \bbZ[\pi_\sigma])=0
\end{equation*} 
when $\sigma$ is an edge or a face, since $\pi_\sigma$ is a 1-ended group ($\pi_\sigma=\bbZ^2$ and $\bbZ^3$). 
\end{proof}

The condition that $\pi$ has tame cohomology is perhaps restrictive (see Remark \ref{tame restrictions}) but we have many examples, including an algorithm that can produce infinitely many \RAAGs\ with tame cohomology. 

\begin{proposition}\label{prop:tame examples}
Let $\Gamma$ be the simplicial graph obtained 
 by carrying out finitely many of the following operations, starting with a single vertex:
\begin{enumerate}
\item Attach a new 1-simplex or a new 2-simplex to the previous graph $\Gamma_0$  by identifying one of its vertices with any (but only one) vertex of $\Gamma_0$.

\item Attach a new  2-simplex to the previous graph $\Gamma_0$ by identifying one of its 1-simplices with any (but only one) 1-simplex of $\Gamma_0$.
\end{enumerate}
Then the \RAAG\ defined by $\Gamma$ will have tame cohomology. 
\end{proposition}

\begin{proof} 
At each step  we add a new simplex to the previously constructed graph $\Gamma_0$,  which can always be
contracted to $\Gamma_0$ in its flag complex. The algorithm  provided   by repeating steps (i) and (ii) produces a  simplicial graph $\Gamma$ with contractible flag complex $\hat\Gamma$, which is clearly $2$-connected. 

 Furthermore, this algorithm guarantees that $H^1(\lk(\sigma))=0$  for the link of every vertex $\sigma$.  We may then apply Lemma \ref{tame} to conclude that the \RAAG\ defined by $\Gamma$ has tame cohomology.
\end{proof}

In the setting of Theorem A, we restrict attention to \RAAGs\ with no $4$-cliques in their defining graphs.
It is possible that all such \RAAGs\ have tame cohomology, but we do not yet know if the first condition is always satisfied.

\begin{example}
Does the \RAAG\ $\pi$ with the following defining graph have $H^2(\pi;\La)^* = 0$~?
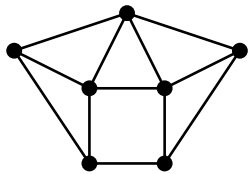
\begin{center}
\begin{tikzpicture}[scale = 1]
\tikzstyle{every node}=[draw, shape=circle, minimum size = 1cm, line width = 1pt];
\path (0,0) coordinate (a);
\path (0,1) coordinate (b);
\path (1,1) coordinate (c);
\path (1,0) coordinate (d);
\path (-1,1.5) coordinate (e);
\path (0.5,2) coordinate (f);
\path (2,1.5) coordinate (g);
\draw[line width = 1pt] (a)--(b)--(c)--(d)--(a)--(e)--(f)--(b)--(e) (f)--(c)--(g)--(f) (g)--(d);
\foreach \i in {a,b,c,d,e,f,g}{\fill (\i) circle (3pt);}
\end{tikzpicture}
\end{center}
\end{example}

 The following examples provide graphs with non-simply connected flag complexes whose  associated \RAAGs\ have  $H^2(\pi;\La)^*=0$.
This shows at least that  Lemma \ref{tame} is not the best possible result.

\begin{example} \label{non-simply connected tame}
The \RAAGs\ defined by the following graphs all have $H^2(\pi;\La)^*=0$. 
\smallskip

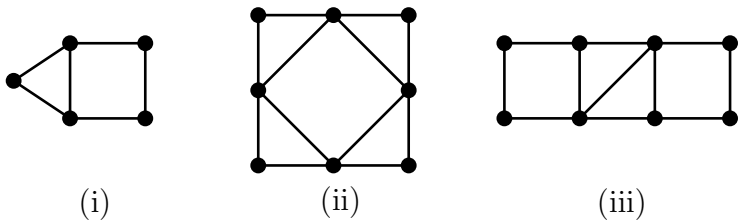
\begin{center} 
\begin{minipage}{.2\textwidth}
\begin{tikzpicture}[scale = 1]
\tikzstyle{every node}=[draw, shape=circle, minimum size = 1cm, line width = 1pt];
\path (0.25,0) coordinate (x1);
\path (1,0.5) coordinate (x2);
\path (1,-0.5) coordinate (x3);
\path (2,0.5) coordinate (x4);
\path (2,-0.5) coordinate (x5);
\path (1,-1.25) coordinate (z); 
\path (1,1) coordinate (p);
\draw[line width = 1pt] (x1)--(x2)--(x3)--(x1) (x2)--(x4)--(x5)--(x3);
\foreach \i in {1,..., 5}{\fill (x\i) circle (3pt);}
\end{tikzpicture} \\
\centering (i)
\end{minipage}
\begin{minipage}{.2\textwidth}
\begin{tikzpicture}[scale = 1]
\tikzstyle{every node}=[draw, shape=circle, minimum size = 1cm, line width = 1pt];
\foreach \i in {1,...,3}{\path (\i,2) coordinate (x\i);}
\foreach \i in {1,3}{\path (\i,1) coordinate (z\i);}
\foreach \i in {1,...,3}{\path (\i,0) coordinate (y\i);}
\draw[line width = 1pt] (x1)--(x3)--(y3)--(y1)--(x1) (x2)--(z3)--(y2)--(z1)--(x2);
\foreach \i in {1,..., 3}{\fill (x\i) circle (3pt); \fill (y\i) circle (3pt);}
\foreach \i in {1,3}{\fill (z\i) circle (3pt);}
\end{tikzpicture}\\
\centering (ii) \quad
\end{minipage}
\begin{minipage}{.25\textwidth}
\begin{tikzpicture}[scale = 1]
\tikzstyle{every node}=[draw, shape=circle, minimum size = 1cm, line width = 1pt];
\foreach \i in {1,...,4}{\path (\i,1) coordinate (x\i);}
\foreach \i in {1,...,4}{\path (\i,0) coordinate (y\i);}
\draw[line width = 1pt] (x1)--(x4)--(y4)--(y1)--(x1) (x2)--(y2)--(x3)--(y3);
\foreach \i in {1,..., 4}{\fill (x\i) circle (3pt); \fill (y\i) circle (3pt);}
\path (1,-0.75) coordinate (z);
\path (1,1.5) coordinate (p);
\end{tikzpicture}\\
\centering (iii) \qquad
\end{minipage}
\end{center}
\end{example}

The argument involves dualizing the short exact sequences derived from the Mayer-Vietoris sequence for $\pi$ expressed as the amalgamated product of two subgroups. We denote by $\pi_i$ the \RAAG\ associated to a graph $\Gamma_i$. We first remark that if $\Gamma$ is the union of two subgraphs $\Gamma_1$ and $\Gamma_2$ which intersect at $\Gamma_0$, then $\pi$ is the amalgamated product of $\pi_1$ and $\pi_2$ over $\pi_0$: this can be seen by comparing the presentations for $\pi_1$ and $\pi_2$ with the relations of $\pi_0$. Secondly, if $\Gamma_i$ is the induced subgraph on a collection of vertices from $\Gamma$, the groups $\pi_i$ are subgroups of $\pi$. 

\begin{proposition} \label{mv breakdown} 
If $\Gamma=\Gamma(\pi)$ satisfies the following conditions, then $H^2(\pi;\La)^*=0$.
\begin{enumerate} \addtolength{\itemsep}{0.3\baselineskip}
\item $\Gamma$ is the union of two graphs $\Gamma_1$ and $\Gamma_2$ which intersect in $\Gamma_0$, and all three subgraphs induce subgroups of $\pi$. 
\item $\Gamma_1$ and $\Gamma_2$ are connected and have 2-connected flag complexes.
\item $\Gamma_0$ is a disjoint union of vertices.
\end{enumerate}
\end{proposition}

\begin{proof}
Since $\Gamma_1$ and $\Gamma_2$ are connected and $\Gamma_0$ is 0-dimensional, we get the short exact sequence
$$ 0 \to H^1(\pi_0;\La) \to H^2(\pi;\La) \to H^2(\pi_1;\La)\oplus H^2(\pi_2;\La) \to 0$$
from Mayer-Vietoris. Since $\hat\Gamma_1$ and $\hat\Gamma_2$ are 2-connected, $H^2(\pi_1;\La)^*$ and $H^2(\pi_2;\La)^*$ are both zero by Lemma \ref{Hom duals are 0}. Additionally, $H^1(\pi_0;\La)^*=0$ by Lemma \ref{lem:disjoint}.
Thus dualizing the short exact sequence above gives the desired result.
\end{proof}

We can apply this proposition to the graph $\Gamma$ in Example \ref{non-simply connected tame} (i) with the following two subgraphs below, to be denoted by $\Gamma_1$ and $\Gamma_2$, respectively:
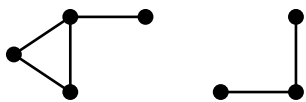
\begin{center}
\begin{tikzpicture}[scale=1]
\tikzstyle{every node}=[draw, shape=circle, minimum size = 1cm, line width = 1pt];
\path (-0.75,0.5) coordinate (x1);
\path (0,0) coordinate (x2);
\path (1,1) coordinate (x4);
\path (0,1) coordinate (x5);
\path (2,0) coordinate (x6);
\path (3,0) coordinate (x7);
\path (3,1) coordinate (x8);
\draw[line width = 1pt] (x1) -- (x2) (x4) -- (x5) -- (x1) (x2) -- (x5) (x6) -- (x7) -- (x8);
\foreach \i in {1,2,4,5,6,7,8}{\fill (x\i) circle (3pt);}
\end{tikzpicture} 
\end{center}
Together, their union $\Gamma_1 \cup \Gamma_2$ is $\Gamma$ and their intersection $\Gamma_0$ is two disjoint vertices. Both $\Gamma_1$ and $\Gamma_2$ are connected and have 2-connected flag complexes. Thus Proposition \ref{mv breakdown} shows the graph in Example \ref{non-simply connected tame} (i) satisfies the first condition for tame cohomology.

\section{Finitely presented groups with $\cd \pi \leq 3$}\label{sec:fpgroups} 

In this section we will show that any finitely presented group $\pi$ with $\cd\pi \leq 3$ satisfies the conditions (ii) and (iii) for tame cohomology. In particular, these conditions hold for a \RAAG\  $\pi$ with no $4$-cliques in its defining graph $\Gamma$. The arguments in this section use the homological algebra definition of the $\Ext$-functors via multiple extensions  (see MacLane \cite[Chap.~III, \S\S 5-6]{MacLane:1967}).

A particularly useful tool is the hypercohomology spectral sequence defined in Benson \cite[Proposition 3.4.3]{benson2}.
If $\bC$ and $\bD$ are chain complexes over a ring $\La$, with $\bC$ bounded above and $\bD$ bounded below, then there is a spectral sequence
$$ \Ext^p_{\La}(H_q(\bC), \bD) \Rightarrow \Ext^{p+q}_{\La}(\bC, \bD)$$
converging to the hypercohomology groups (see Benson and Carlson \cite[\S 2.7]{benson1} for the definition of $\Ext$ for chain complexes). This is a first quadrant spectral sequence, so the indexing convention above is standard ($p$ runs along the $x$-axis). The differential $d_r$ has bi-degree $(r, -r+1)$. 

A summary of the information we need about this spectral sequence is provided for the reader's convenience in Section \ref{sec:append}. 

\begin{proposition} \label{almost tame} 
Let $\pi$ be a finitely presented group with $\cd \pi \leq 3$.
Then $H^3(\pi;\La)^*=0$ and $\Ext^1_\La(H^3(\pi;\La), \La)=0$. 
\end{proposition}

\begin{proof}
By Wall \cite[Theorem E]{wall-finiteness1}, a finitely presented group $\pi$ with $\cd \pi \leq 3$ has geometric dimension $\leq 3$. In other words, there exists a $3$-dimensional, finite, aspherical complex $L$ with $\pi_1(L) = \pi$. Let $Y=Y(L)$ be an 8-dimensional thickening of $L$, and let $Z$ be the double of $Y$. The long exact sequence of the pair $(Y,\partial Y)$ with $\La$-coefficients gives isomorphisms $H_{i-1}(\partial Y) \to H_{i-1}(Y)$ for $i\leq 4$, since $H_i(Y,\partial Y) \cong H^{8-i}(Y) \cong H^{8-i}(L)$. Furthermore, $H_{i-1}(Y) = 0$ for $i>1$. Then by the Mayer-Vietoris sequence, $H_i(Z;\La)=0$ for $1<i \leq 4$, and $H_1(Z;\La)=0$ since $\widetilde Z$ is simply connected.

We use the \textbf{universal coefficient spectral sequence}  which has $E_2$-page 
$$E^{p,q}_2 = \Ext^p_\La(H_q(Z;\La),\La).$$ 
This is a special case of the hypercohomology spectral sequence, where $\bC = C(Z;\La)$ and $\bD = \La$ (a $0$-dim chain complex with $\La$ in degree zero). 

Since, $H_q(Z;\La)=0$ for $q=1,\dotsc,4$, the $E^{p,q}_2$-terms are all zero for $1\leq q \leq 4$, and since $H^p(\pi;\La)=0$ for all $p>3$, the $E^{p,q}_2$-terms are zero for $p>3$. Thus no differentials affect the $E^{0,5}_2$- and $E^{1,5}_2$-terms, which are $E^{0,5}_2 \cong H^3(\pi;\La)^*$  and $E^{1,5}_2 \cong \Ext^1_\La(H^3(\pi;\La),\La)$ via the isomorphisms $H_5(Z;\La) \cong H^3(Z;\La) \cong H^3(\pi;\La)$. Since $H^5(Z;\La) \cong H_3(Z;\La) = 0$ and $H^6(Z;\La) \cong H_2(Z;\La) = 0$, these two terms must be zero. 
\end{proof}

In the universal coefficient spectral sequence computing $H^*(Z;\La)$, since $H^7(Z;\La) = H_1(Z;\La) = 0$,   the differential 
\eqncount
\begin{equation}\label{eq:differential}
d_2\colon E^{0,6}_2 = H^2(\pi;\La)^* \xrightarrow{\cong} E^{2,5}_2 = \Ext^2_\La(H^3(\pi;\La),\La)
\end{equation}
is an isomorphism.
In addition,  $\Ext^1_\La(H^2(\pi;\La))$ injects into $\Ext^3_\La(H^3(\pi;\La),\La)$. 
 
From the universal coefficient spectral sequence for the $2$-skeleton $K \subset L$, we get the exact sequence:
\eqncount
\begin{equation}\label{seq for K}
0 \to H^2(\pi;\La) \to H^2(K;\La) \to \pi_2(K)^* \to H^3(\pi;\La) \to 0
\end{equation}
where
the last non-zero map is a surjection since $H^3(K;\La)=0$.
This sequence can be analyzed 
by splicing together the short exact sequences:
$$0 \to H^2(\pi;\La) \to H^2(K;\La) \to V \to 0
\quad  \text{and} \quad
0 \to V \to \pi_2(K)^* \to H^3(\pi;\La) \to 0,$$
where 
$$V:= \cok(H^2(\pi;\La) \to H^2(K;\La)) = \ker(\pi_2(K)^* \to H^3(\pi;\La)).$$
After dualizing,  the associated long exact cohomology sequences give
 two connecting homomorphisms $\delta_1\colon H^2(\pi;\La)^* \to \Ext_\La^1(V, \La)$ and
$\delta_2\colon \Ext_\La^1(V, \La) \to \Ext_\La^2(H^3(\pi;\La), \La)$. 

 \begin{lemma}\label{lem:differential}
Let $\pi$ be a finitely presented group with $\cd \pi \leq 3$. 
The differential  \textup{\eqref{eq:differential}} induces a natural isomorphism
$d_2\colon H^2(\pi;\La)^* \cong \Ext^2_\La(H^3(\pi;\La),\La)$,
where $d_2 = \delta_2\circ \delta_1$.
\end{lemma}
\begin{proof} We refer to Section \ref{sec:append} for the $d_2$ differential formula \eqref{one} given by  splicing, and apply it to the manifold
 $Z = Y(L)\, \cup\, Y(L)$ defined above. Note that $\bC= C_*(Z;\La) \simeq C^{8-*}(Z;\La)$ by Poincar\' e duality. We claim that the sequence \eqref{one}  for $p=0$ and $q=6$, and expressed in terms of the cochain complex of $Z$, is exactly the sequence \eqref{seq for K}. 

Here are some additional details to explain this identification.
Let $\bD := C^{8-*}(Z;\La)$ denote the cochain complex, and  compute the same sequence \eqref{one} for the complex $\bD$ (with $q=6$).
Then $H_6(\bD) = H^2(Z;\La)$ and $H_5(\bD) = H^3(Z;\La)$. It is not hard to see (by definition of cohomology for $Z$) that the term 
$D_6/B_7 \cong H^2(K;\La)$. In fact,   the sequence \eqref{one} in this case is just
$$ 0 \to H^2(Z;\La) \to H^2(K;\La) \to H^3(Z,K;\La) \to H^3(Z;\La) \to 0,$$
which is isomorphic to the sequence \eqref{seq for K} 
since $H^3(Z,K;\La) \cong \pi_2(K)^*$, by an application of the universal coefficient spectral sequence to the pair $(Z,K)$ (the only term on the 3 line is at (0,3)).

  Since $C \simeq D$, the corresponding   $d_2$ in the spectral sequence for computing $H^*(D)$ is given by splicing with  \eqref{seq for K} as claimed.
 Finally, note that the formula $d_2 = \delta_2\circ \delta_1$ for $d_2\colon H^2(\pi;\La)^* \cong \Ext^2_\La(H^3(\pi;\La),\La)$ is obtained by composing the connecting maps $\delta_1$ and $\delta_2$ obtained from two 3-term exact sequences which splice together to give the 4-term exact sequence  \eqref{seq for K}.
\end{proof}

\begin{remark}
These results, which apply to all \RAAGs\ with $H_4(\pi;\bZ)=0$, are sharper in this special case than our conclusions about tame cohomology from Section \ref{sec:tame}. For example, we no longer require the associated flag complex $\hat\Gamma$ to be simply connected.
\end{remark}

We are indebted to one of the referees for pointing out that Proposition \ref{almost tame}
 is a special case of a more general cohomological result, with an analogous proof.  Here is the statement together with an alternate proof.

\begin{proposition}
Let $\pi$ be a finitely presented group with $\cd \pi \leq n$.
Then $H^n(\pi;\La)^*=0$ and $\Ext^1_\La(H^n(\pi;\La), \La)=0$. 
\end{proposition} 

\begin{proof}
We begin with the hypercohomology spectral sequence \cite[Proposition 3.4.3]{benson2}, which has $E_2$-page
$$E_2^{p,q} = \Ext_\La^p(H^{-q}(\pi;\La),\La).$$
This is a fourth quadrant spectral sequence which converges to $H_{p+q}(\pi;\La)$. 
For all group homology with group ring coefficients, $H_{p+q}(\pi;\La)=0$ if $p+q\neq 0$ and $\bbZ$ if $p+q=0$, and so all entries  not lying on the line $p+q=0$ will die out in the spectral sequence. However, since $\pi$ has $\cd \pi \leq n$, the terms $E_2^{p,q}=0$ when $q\leq -(n+1)$. There are no differentials that affect $E^{0,-n}_2=H^n(\pi;\La)^*$ and $E^{1,-n}_2=\Ext^1_\La(H^n(\pi;\La),\La)$, and so both must be zero.
\end{proof}

\begin{remark}
The hypercohomology spectral sequence above, combined with Lemma \ref{lem:disjoint}, shows that $H^1(\pi;\La)^* \neq 0$ for a \RAAG\ $\pi$ if and only if its defining graph is disconnected and contains at least one edge. The simplest example is $\pi = \bZ^2 \ast \bZ$. For the groups $\pi=\pi_{n.m}$, we let $Z$ denote the double of an 8-dimensional thickening of the Salvetti complex. Since $H^8(Z;\La) = \bZ$, we obtain an exact sequence
$$0 \to H^1(\pi;\La)^* \to \Ext_\La^2(H^2(\pi;\La), \La) \to \bZ \to \Ext_\La^1(H^1(\pi;\La), \La) \to 0.$$
A further calculation with the filtration sequences shows that 
$\Ext_\La^2(H^2(\pi;\La), \La)$ is a direct sum of $n$ copies of $\Ext_\La^2(\bZ[\pi/\pi_\sigma], \La)$, where $\sigma$ is an edge. If $n>0$ this is enough to show that $H^1(\pi;\La)^* \neq 0$. If  $m=0$ as well, then   the exact sequence above reduces to 
$$0\to H^1(\pi;\La)^* \to (\bZ[\pi/\pi_\sigma])^n \to \bZ \to 0, $$
 and $H^1(\pi;\La)^* \cong \La^{n-1}$.  For each of the $n$ edges,  $\pi_\sigma \cong \bZ^2$.
\end{remark}

\section{Minimal models and stabilization} \label{sec:four}
For any finitely presented group $\pi$, one can construct a 4-manifold $M$ with fundamental group $\pi$ by doubling a thickening of a finite 2-complex $K$ with $\pi_1(K)=\pi$. If $K$ has minimal Euler characteristic, then the double of $K$ will have minimal Euler characteristic over all double constructions.

If $\pi$ is a \RAAG, we can take $K$ to be the 2-skeleton of the Salvetti complex mentioned in Section \ref{RAAG section}. The Salvetti complex has minimal Euler characteristic over all possible $\Bpi$ since its chain complex gives a minimal resolution for $\pi$.

\begin{definition}
We say that a $4$-manifold  $X$  is \emph{minimal for} $\pi$   if its Euler characteristic  is minimal over all closed, oriented 4-manifolds with fundamental group $\pi$. 
\end{definition}
The Euler characteristic of a minimal $4$-manifold for $\pi$ is  the Hausmann-Weinberger invariant \cite{Hausmann:1985}. It has been determined for free abelian groups by Kirk and Livingston
\cite{Kirk:2005}, but is still unknown for most classes of finitely presented groups.
\begin{theorem}[Hildum {\cite[Theorem 1.2]{Hildum:2015}}] 
 Let $\pi$ be a \RAAG\ with $H_4(\pi;\bZ) = 0$. 
Let $K$ be the 2-skeleton of the Salvetti complex with fundamental group $\pi$, and let $Y(K)$ denote a spin $4$-dimensional thickening of $K$. Then the double
$M_0 := Y(K)\cup Y(K)$ is a minimal spin $4$-manifold for $\pi$.
\end{theorem}
 \begin{corollary} 
Let $\pi$ be a \RAAG\  with $H_4(\pi;\bbZ)=0$. 
Then 
\begin{equation*}
\rk_\bbZ(\pi_2(M_0) \otimes_\La \bbZ) = b_2(\pi) + b_3(\pi).
\end{equation*}
if $M_0$ is a minimal spin $4$-manifold for $\pi$.
\end{corollary}
\begin{proof}
The  spectral sequence (of the covering) converging to $H_*(M_0)$ yields the following exact sequence:
\begin{equation*}
\xymatrix{0 \ar[r] & H_3(\pi) \ar[r]^{d^3 \quad} \ar[d]^\cong & \pi_2(M_0) \otimes_\La \bbZ \ar[r] \ar[d]^\cong & H_2(M_0) \ar[r] \ar[d]^\cong & H_2(\pi) \ar[r] \ar[d]^\cong & 0 \\ 0 \ar[r] & \bbZ^{b_3(\pi)} \ar[r] & \pi_2(M) \otimes_\La \bbZ \ar[r] & \bbZ^{2b_2(\pi)} \ar[r] & \bbZ^{b_2(\pi)} \ar[r] & 0}
\end{equation*}

The second term in the exact sequence above is the $(0,2)$ term in the spectral sequence, $H_0(\pi;\pi_2(M_0)) \cong \pi_2(M_0) \otimes_\La \bbZ$. By exactness, the $\bbZ$-rank of $\pi_2(M_0) \otimes_\La \bbZ$ is $b_2(\pi) + b_3(\pi)$. 
\end{proof}
   
  The thickened double construction will allow us to determine the structure of $\pi_2(M)$ as a $\La$-module, for any $4$-manifold with fundamental group $\pi$ of cohomological dimension $\cd \pi \leq 3$. First we discuss  the homology and cohomology of  finite $2$-complexes. 
  
  \begin{lemma} \label{H_2(K) stably free} Let $\pi$ be a finitely presented group with $\cd \pi \leq 3$.
If $K$ is any finite $2$-complex
with fundamental group $\pi$, then $\pi_2(K) = H_2(K;\La)$ is stably free as 
 a  $\La$-module.   
\end{lemma}
\begin{proof} 
Let $X$ be an aspherical $3$-complex model for $K(\pi,1)$,
and let $K$ be the 2-skeleton of $X$. In this special case we will show that $H_2(K;\La)$ is a finitely generated \emph{free} $\La$-module. 
Consider the  chain complex $C_*(X):= C_*(X;\La) $ with  $\La$-module coefficients:
\eqncount
\begin{equation}\label{eq:salvetti}
\xymatrix{
0 \ar[r] & C_3(X) \ar[r] & C_2(X) \ar[r] & C_1(X) \ar[r] & C_0(X) \ar[r] & \bbZ \ar[r] & 0 }
\end{equation}
This  chain complex is exact  since $\wX$ is contractible,  and so $H_2(K;\La) \cong Z_2(K) = Z_2(X) = C_3(X)$. Since $C_i(X) \cong \La^{b_i(\pi)}$, we have $H_2(K;\La) \cong \La^{b_3(\pi)}$, where $b_i(\pi)$ denotes the number of $i$-cells in $X$. The result now follows, since the module $H_2(K;\La)$ is independent of the choice of finite $2$-complex $K$, after stabilization with free modules.  
\end{proof}

In general, the available information about $H^2(K;\La)$ is contained in the 4-term exact sequence (\ref{seq for K}):  
\begin{equation*} 
0 \to H^2(\pi; \La) \to H^2(K; \La) \to \Hom_\La(H_2(K;\La),\La) \to H^3(\pi;\La) \to 0.
\end{equation*}
 In certain cases we can make specific calculations.  
\begin{example}\label{ex:notfree}
Consider the case of the \RAAG\ $\pi = \bbZ^3$ in which the associated graph $\Gamma$ is a 3-clique and the flag complex $\hat\Gamma$ is a single 2-simplex. Then, 
 $H^2(\pi;\La)=0$ and $H^3(\pi;\La)=\bbZ$ by Theorem \ref{group ring coefficients}. This example indicates that even though $H_2(K;\La)$ is a free $\La$-module, $H^2(K;\La)$ may not be free. In this case, (\ref{seq for K}) becomes 
\begin{equation*}
0 \to H^2(K;\La) \to \La \to \bbZ \to 0,
\end{equation*}
which shows that $H^2(K;\La)$ is the augmentation ideal $I(\pi)$.
 Thus, $\pi_2(M_0)= I(\pi) \oplus \La$   for $\pi=\bbZ^3$, and 
the module $I(\pi)$ is not (stably) free since $H_1(\pi;I(\pi))\cong H_2(\pi;\bbZ) \cong \bbZ^3$.
\end{example}

We now consider the structure of $\pi_2(M)$ and $s_M$ for the thickened doubles of finite $2$-complexes (see Kreck and Schafer \cite[\S II]{kreck-schafer1}).

\begin{lemma}\label{lem:double} Let $\pi$ be a finitely-presented group and let $K$ be any  finite $2$-complex with fundamental group $\pi$. If $Y(K)$ denotes a  spin  $4$-dimensional thickening of $K$, 
 then 
 $$\pi_2(M) = H^2(K;\La) \oplus H_2(K; \La)$$
  as a $\La$-module, for the thickened double $M = Y(K)\cup Y(K)$.  The equivariant intersection form $s_M$ is a metabolic form. 
\end{lemma}
\begin{proof}
Let $Y=Y(K)$ and notice that $M = \bd(Y \times I)$. We start with the long exact sequence in homology for the pair $(M,Y)$:
\begin{equation*}
\dotsm \to H_3(M,Y;\La) \to H_2(Y;\La) \to H_2(M;\La) \to H_2(M,Y;\La) \to H_1(Y;\La) \to \dotsm
\end{equation*}
The inclusion of $Y$ into $M$ induces a split injective map $H_i(Y;\La) \to H_i(M;\La)$ in every dimension. Thus the maps from $H_i(M,Y;\La)$ to $H_i(Y;\La)$ are all zero maps. In addition,
using excision properties as well as Poincar\'e-Lefschetz duality, we have the isomorphisms $H_2(M,Y;\La) \cong H^2(Y,\partial Y;\La) \cong H_{2}(Y;\La)$.
 This gives the split short exact sequence
 \eqncount
\begin{equation}\label{seq:splitshort}
0 \to H^2(K;\La) \to H_2(M;\La) \to H_2(K;\La) \to 0.
\end{equation}
This, along with the Hurewicz isomorphism $H_2(M;\La) \cong H_2(\wM;\bbZ) \cong \pi_2(M)$ yields the desired result. By \cite[Proposition II.4]{kreck-schafer1}, the equivariant intersection form has the structure:
\eqncount
\begin{equation}\label{eq:doubleform}
s_M((u, x), (v, y)) = u(y) + \overline{v(x)} + \delta_K(u,v), \end{equation}
for all $ u,v \in H^2(K;\La)$ and $x,y \in H_2(K;\La)$. In this formula, $\delta_K$
 is an even hermitian symmetric form on $H_2(K;\La)$ determined by the thickening. 
\end{proof}

We can now determine the structure of $\pi_2(M)$ in general, after stabilization by free $\Lambda$-modules.
 \begin{proposition}
\label{prop:stable}
 Let $N$ be a closed, oriented, spin, TOP $4$-manifold 
 with 
 fundamental group $\pi$. If $\cd \pi\leq 3$, then there exists a simply connected, closed $4$-manifold $X$ such that 
$$N' \Sharp r(S^2 \times S^2) \approx  M' \Sharp s(S^2 \times S^2),$$
 for some $r, s\geq 0$, where
$N':=N \Sharp \CP^2 \Sharp \overline{\CP^2}$,  $M':=M_0 \Sharp \CP^2 \Sharp \overline{\CP^2}\Sharp X$, and
$M_0$ is a thickened double.
\end{proposition}

\begin{proof}
We note that the thickened   double model $M_0$ represents the zero bordism element in  $\Omega_4^{STOP}(K(\pi,1)) = \bZ$. Since the signature detects elements in this bordism group, it follows that $N':=N \Sharp \CP^2 \Sharp \overline{\CP^2}$ is stably homeomorphic to $M':=M_0 \Sharp \CP^2 \Sharp \overline{\CP^2}\Sharp X$ for some closed, simply connected $4$-manifold $X$ with $\sign(N') = \sign(X)$. If  $M_0$ is minimal we can take $X$ to be a connected sum of copies of $\CP^2$  and $\overline{\CP^2}$  so that $N'$ and $M'$ have the same signature and Euler characteristic. It follows that
$N' \Sharp r(S^2 \times S^2) \approx  M' \Sharp r(S^2 \times S^2)$, for some $r \geq 0$ (see  \cite[Theorem C]{Kreck99} and  \cite[\S 9.1]{Freedman:1990}).
\end{proof}

\begin{corollary}\label{cor:stableiso}
Let $M$ and $N$ be a closed, oriented, spin, TOP $4$-manifolds with
fundamental group $\pi$. If $\cd \pi \leq 3$, then 
$\pi_2(N)$ is stably isomorphic to $\pi_2(M)$ as a $\La$-module.
\end{corollary}
\begin{proof}
By Proposition \ref{prop:stable}, the modules $\pi_2(N)$ and $\pi_2(M)$ are both stably isomorphic to $\pi_2(M_0)$ after stabilization by direct sum with free $\La$-modules.  
\end{proof}

\begin{remark}
By using thickenings, we have shown that $\pi_2(N)$ is stably isomorphic to $H_2(K;\Lambda) \oplus H^2(K;\Lambda)$, for any $K$ any finite $2$-complex with $\pi_1(K) = \pi_1(N)$. This is a special case of a general result (see \cite[Proposition 2.4]{HambletonKreck88},  \cite[Theorem 4.2] {h5}), that $\pi_2(N)$ is stably isomorphic to a $\Lambda$-module $E$ in an extension:
$$ \cE_N : 0 \to H_2(K;\Lambda) \to E \to  H^2(K;\Lambda) \to 0\ .$$
For fundamental groups with $\cd \pi \leq 3$, our direct argument verifies that $\cE_N$ is stably a split extension.
\end{remark}

\begin{remark} \label{rem:k-invariant}
The $k$-invariant of the thickened   double $M = Y(K) \cup Y(K)$ is the image (induced by the inclusion) of the $k$-invariant of $K$.  More generally, if $i_M\colon K\subset M$ is the $2$-skeleton of a closed topological $4$-manifold $M$, then $(i_M)_*(k_K) = k_M \in H^3(\pi; \pi_2(M))$. This follows from the push-out diagram:
$$\xymatrix{
0 \ar[r] &\pi_2(K) \ar[r] \ar[d]^{(i_M)_*}& C_2(M) \ar[r] \ar[d]& C_1(M) \ar[r] \ar@{=}[d]& C_0(M) \ar[r] \ar@{=}[d]& \bbZ \ar[r] \ar@{=}[d]& 0\cr
0 \ar[r] & \pi_2(M) \ar[r] & C_2(M)/B_2(M) \ar[r] & C_1(M) \ar[r] & C_0(M) \ar[r] & \bbZ \ar[r] & 0
}$$
in which the upper multiple extension represents $k_K$ and the lower extension represents $k_M$ (see Baues \cite[p.~101]{Baues:1991}).
\end{remark}

Here is a useful consequence of this description for the $k$-invariant.
\begin{proposition}\label{prop:movek}
 Let  $M$ and $N$ be  closed topological $4$-manifolds with fundamental group $\pi$ of $\cd \pi \leq 3$. If the equivariant intersection forms $s_M$ and $s_N$ are stably isometric, then their quadratic $2$-types $Q(M)$ and $Q(N)$ are stably isometric.
\end{proposition}

\begin{proof}    
It is enough to show that  we have an isometry $Q(M) \cong Q(N)$ after stabilization by connected sum with sufficiently many copies of $S^2 \times S^2$.  After preliminary stabilization, we may assume that $s_M \cong s_N$ and that $i_M\colon K \subset M$ and $i_N\colon K \subset N$ is a $2$-skeleton for both manifolds. This uses Whitehead's stable uniqueness theorem for finite $2$-complexes with a fixed fundamental group (up to simple homotopy equivalence) after wedging with sufficiently many copies of $S^2$.  By Remark \ref{rem:k-invariant}, we have $(i_M)_*(k_K) = k_M $ and
$(i_N)_*(k_K) = k_N $. By general position, we may assume that the $1$-skeleton of $K$ is embedded in $M$ and the $2$-cells are all immersed in $M$,  with at most  double point intersections and self-intersections.

Since $\pi$ has $\cd \pi \leq 3$, we may assume by Lemma \ref{H_2(K) stably free} that $\pi_2(K)$ is a finitely generated free $\La$-module, and a direct summand of both $\pi_2(M)$ and $\pi_2(N)$. Finally, after further stabilization, we may assume that the equivariant intersection forms $s_M$ and $s_N$ and their quadratic refinements $q_M$ and $q_N$ vanish on these summands $\pi_2(K)$.
This is achieved as follows: at each $2$-cell intersection or self-intersection point 
 where two sheets meet, we take the connected sum with a copy of $S^2 \times S^2$. Then we construct a new embedding of $K$ by attaching disjoint $2$-cells in $S^2 \times S^2 \setminus D^4$ to the original separate $2$-cell sheets. This process removes one intersection point, and may be repeated to embed the same complex $K$ in a stabilization of $M$. 
 
 After applying the same steps to the image of $K$ in $N$, we may assume that
$\pi_2(K)$ is a \emph{totally isotropic} submodule with respect to both quadratic forms. At this stage, $K$ is no longer the full $2$-skeleton of $M$ and $N$, but still carries their $k$-invariants.  

Now suppose that $\phi\colon (\pi_2(M), s_M) \to (\pi_2(N), s_N)$ is an isometry of  equivariant intersection forms. We claim that there is a self-isometry $\psi \colon s_N \to s_N$, after further stabilization,  such that the composite $\psi\circ\phi$ satisfies $\psi(\phi(x))= x$, for all $x \in \pi_2(K) \subset \pi_2(M)$. It will then follow that $(\psi\circ \phi)_*(k_M) = k_N$ and hence that $Q(M) \cong Q(N)$. 

To construct the required self-isometry of a stabilized $s_N$, we use unitary transvections as in \cite[\S 1]{hk5}. 
Let $\{w_1, \dots, w_k\}$ be a free $\Lambda$-base for $\pi_2(K)$, which can be considered either as elements in $\pi_2(M)$ or in $\pi_2(N)$. Since $\pi_2(N)^* \to \pi_2(K)^*$ is a split surjection, we can pick elements $\{v_1, \dots, v_k\}\subset \pi_2(N)$ such that $\{\ad s_N(v_1), \dots, \ad s_N(v_k)\}$ projects to the dual base $\{w_1^*, \dots, w_k^*\}$ for  $\pi_2(K)^*$. In other words, $s_N(w_i, v_j) = \delta_{ij}$, for $1\leq, i, j \leq k$. Since $s_N$ admits a quadratic refinement $q_N$,  and $\pi_2(K)$ is a free totally isotropic $\La$-module, we may assume that $q_N(v_i) = 0$ and $s_N(v_i, v_j) = 0$, for $1\leq i, j \leq k$ (see \cite[Lemma 5.3]{wall-book}). It follows that the elements $\{w_1, \dots,w_k; v_1, \dots, v_k\}$ span a hyperbolic direct summand  of $(\pi_2(N), s_N, q_N)$ isomorphic to $H(\pi_2(K))$. 

We consider the orthogonal direct sum 
$s_N \oplus  H(\Lambda)_1 \oplus \dots \oplus H(\Lambda)_k$
with $k$ new hyperbolic summands,  and let $\{e_i,f_i\}$ denote a hyperbolic base for the $i$-th summand.
 We compute   the following composite of unitary transvections (see the formula  \cite[1.2]{hk5}):
$$\theta_i(w_i) := (\sigma_{f_i, 0,w_i}\circ \sigma_{e_i,0,v_i})(w_i) = \sigma_{f_i, 0,w_i}(\sigma_{e_i,0,v_i}(w_i)) =\sigma_{f_i,0,w_i}(w_i+e_i) = e_i$$
for $1 \leq i \leq k$.
This shows that   we can move any basis element $w_i$ to the standard position $e_i$   in a new hyperbolic summand. By construction, $\theta_i(w_j) = w_j$, $\theta_i(e_j) = e_j$ and $\theta_i(f_j) = f_j$, if $j \neq i$. Let $\theta = \theta_k \circ \dots \circ \theta_1$, so that $\theta(w_i) = e_i$, for $1\leq i \leq k$.

We can apply this process to elements $\{\phi(w_1), \dots, \phi(w_k)\}$ of the totally isotropic submodule $\phi(\pi_2(K)) \subset \pi_2(N)$, and construct a self-isometry $\theta'$ so that $\theta'(\phi(w_i)) = e_i$, for $1\leq i \leq k$. 
The  required stabilized self-isometry $\psi$ of $s_N$, such that $\psi(\phi(x)) = x$ for all $x \in \pi_2(K)$, is given by $\psi:= \theta^{-1}\circ \theta'$. \end{proof}


\section{The reduced equivariant intersection form}\label{sec:five}

In this section, let $\pi$ denote a finitely presented group with $\cd \pi \leq 3$. 
Let $M$ denote a closed, oriented, spin 4-manifold with fundamental group $\pi$.
As above, we can construct a model $4$-manifold $M_0$ by taking the double of a thickening $Y(K)$, where $K$ denotes the $2$-skeleton of a minimal aspherical $3$-complex $L$ with $\pi_1(L) = \pi$.

For fundamental groups of geometric dimension $\leq 2$, the quotient  
$$\pi_2(M)^\dagger := \pi_2(M) / H^2(\pi;\La)\cong \Hom_\La(H_2(M;\La),\La)$$
 is stably free  since $H^3(\pi;\La)=0$ (see \cite[Corollary 4.4]{HKT09}). 
 In that context, the  ``reduced" equivariant intersection form ${s_M}^\dagger$ was a non-singular form on $\pi_2(M)^\dagger$. 
In our setting, where $\pi$ can have geometric dimension $3$, the quotient $\pi_2(M)^\dagger$  is $\La$-finitely generated   but may not be stably free as a $\La$-module. Instead, we define 
$$\LM:= \pi_2(M)^* = \Hom_\La(H_2(M;\La),\La), $$
 and obtain a non-singular hermitian form  $h_M\colon\LM \times \LM \to \La$, which we also call the \emph{reduced equivariant intersection form}. For groups with $\cd \pi \leq 3$, this form has a key role in the proof of Theorem \ref{thm:maintame}, through a generalization of \cite[Theorem 5.13]{HKT09}.

\begin{lemma}\label{s_M star}
Let $\pi_1(M)$ be a finitely presented group with $\cd\pi_1(M)\leq 3$.
 Then 
\begin{enumerate}
\item $(\ad s_M)^*$ is an isomorphism of finitely generated, stably free  $\La$-modules, and 
\item the inverse of $(\ad s_M)^*$ is the adjoint of a non-singular
 form
\begin{equation*}
h_M\colon  \LM \times \LM  \to \La,
\end{equation*} 
 on $\LM =\pi_2(M)^*$.
\end{enumerate}  
\end{lemma}
\begin{proof} 
 By Corollary \ref{cor:stableiso}, we may assume that $M=M_0$ is a model manifold with the given fundamental group, obtained as the double of a thickened finite $2$-complex $K$.  We will first show that $L_M$ is a finitely generated, (stably) free  $\La$-module after dualizing the summands in the isomorphism 
  $$\pi_2(M)\cong \pi_2(K) \oplus H^2(K;\La)$$
  provided by Lemma \ref{lem:double}. 
  Consider the 4-term exact sequence (\ref{seq for K}) for the $2$-skeleton  $K$ of an aspherical $3$-complex:
\begin{equation*}
0 \to H^2(\pi; \La) \to H^2(K; \La) \to \Hom_\La(H_2(K;\La),\La) \to H^3(\pi;\La) \to 0.
\end{equation*}
The third non-zero term is a free $\La$-module, since $H_2(K;\La)\cong \La^{b_3(\pi)}$. For the following argument, we define $F:=\Hom_\La(H_2(K;\La),\La) \cong \La^{b_3(\pi)}$. 
Recall from Section \ref{sec:fpgroups} that the above sequence is obtained by splicing together 
$$0\to H^2(\pi;\La) \to H^2(K;\La) \to V \to 0$$
and
$$0 \to V \to F \to H^3(\pi;\La) \to 0,$$
where $V:= \cok(H^2(\pi;\La) \to H^2(K;\La)) = \ker(F \to H^3(\pi;\La))$.
Taking the dual of both short exact sequences yields 
$$0 \to V^* \to H^2(K;\La)^* \to H^2(\pi;\La)^* \xrightarrow{\delta_1} \Ext_\La^1(V, \La) \to \dotsm$$
and 
$$0 \to H^3(\pi;\La)^* \to F^* \to V^* \xrightarrow{} \Ext_{\La}^1(H^3(\pi;\La),\La) \to \dotsm\ .$$
By Proposition \ref{almost tame}, we have $ F^* \cong V^* $, and since $F$ is a free $\La$-module the connecting map
$\delta_2\colon \Ext_\La^1(V, \La) \to \Ext_\La^2(H^3(\pi;\La), \La)$
is an isomorphism. 
Now Lemma \ref{lem:differential} shows that the connecting map $\delta_1\colon H^2(\pi;\La)^* \to \Ext_\La^1(V, \La)$   is also an isomorphism. It follows that $V^* \cong H^2(K;\La)^*$.
Together these facts give an isomorphism $H^2(K;\La)^*\cong F^* $ to  a free $\La$-module. 
We have now shown that 
$$\pi_2(M)^* \cong F \oplus F^* \cong \La^{2b_3(\pi)}$$
is a finitely generated, free $\La$-module.

  To show that $(\ad s_M)^*$ is an isomorphism, we can perform the same splicing trick for the 4-term exact sequence 
\eqncount
\begin{equation}\label{seq for M}
0 \to H^2(\pi; \La) \to \pi_2(M) \xrightarrow{\ad s_M} \pi_2(M)^* \to H^3(\pi;\La) \to 0
\end{equation}
with $V_1$ denoting the similar cokernel
$$0 \to H^2(\pi;\La) \to H^2(M;\La) \to V_1 \to 0.$$
We may compare the $4$-term exact sequences for $H^2(M;\La)$ and $H^2(K;\La)$, and note that the restriction map $i^*\colon H^2(M;\La)\to H^2(K;\La)$ is a split surjection by Lemma \ref{lem:double}.
The connecting maps $\delta_1\colon H^2(\pi;\La)^* \to \Ext^1(V;\La)$ and $\delta'_1\colon H^2(\pi;\La)^* \to \Ext^1(V_1;\La)$ are related by the formula $\delta'_1 = i^*\circ \delta_1$. Since the induced map 
$i^*\colon \Ext^1_{\La}(V;\La)\to \Ext^1_{\La}(V_1;\La)$ is a (split) injection and $\delta_1$ is an isomorphism (by Lemma \ref{lem:differential}), we conclude that $\delta'_1$ is an injection. From the sequence
$$0 \to V_1^* \to H^2(M;\La)^* \to H^2(\pi;\La)^* \xrightarrow{\delta'_1} \Ext_\La^1(V_1, \La) \to \dotsm$$
we obtain 
$$V_1^* \cong H^2(M;\La)^* = \pi_2(M)^* = L_M$$
This 
yields the isomorphism 
\eqncount
\begin{equation}\label{seq for ad sM}
(\ad s_M)^* \colon L_M^* \cong (\pi_2(M)^{*})^* \xrightarrow{\cong} V_1^* \xrightarrow{\cong} \pi_2(M)^* = L_M.
\end{equation}
Thus we can define a non-singular
  form $h_M\colon \LM\times \LM\to \La$ by the formula
$$h_M(x,y) = ((\ad s_M)^*)^{-1}(x) (y) \in \La,$$
whose adjoint is the inverse of $(\ad s_M)^*\colon (\pi_2(M)^*)^* \to \pi_2(M)^*$. 
\end{proof}

It remains to check that the form $h_M$ is hermitian.
\begin{definition} Let $h \colon B \times B \to \La$ be a sesquilinear form on a right $\La$-module $B$. We say that $(B,h)$ is \emph{hermitian} (symmetric) if $h = Th$, where  $Th(b,b') :=\overline{h(b', b)}$, for all $b, b'\in B$. 
We say that $h$ is \emph{strongly even}  if $h = \lambda + T\lambda$ for some sesquilinear form $(B, \lambda)$. 
\end{definition}

\begin{remark}
It is a standard fact that an even hermitian form $(B,h)$ is strongly even whenever $B$ is a finitely generated free $\La$-module. 
\end{remark}

Let $\gamma \colon B \to (B^*)^*$ denote the map  $\gamma(b) = b^{**}$ defined by the evaluation
$\la\gamma(b), v\ra = \la b^{**}, v\ra = \overline{v(b)}, \ \text{for\ } v \in B^*$
and note that we have a map of right $\La$-modules:
$$ \langle \gamma(b\lambda),v\ra := \gamma(b\lambda)(v) = \overline{v(b\lambda)} = \overline{(v\bar\lambda)(b)} = \gamma(b)(v\bar\lambda) = (\gamma(b)\cdot \lambda)(v), $$ 
where as usual, the right $\La$-module structure on $B^* = \Hom_\La(B, \La)$ is given by the formula $(v\lambda)(b) = v(b\overline{\lambda})$, for all $v \in B^*$ and $\lambda \in \La$.

  \begin{lemma}\label{lem:hermitian}
   The form $(B,h)$ is hermitian if and only if $(\ad h)^* \circ \gamma = \ad h$.
\end{lemma}
\begin{proof} It suffices to show that $\ad (Th) = (\ad h)^* \circ \gamma$. This follows from the calculation 
\eqncount
\begin{equation}\label{eq:gamma1}
 \la ((\ad h)^* \circ \gamma) (b), b'\ra = \la \gamma(b), \ad h(b')\ra = \overline{\ad h(b')(b)} = \overline{h(b',b)}, 
\end{equation}
where $\la -,-\ra\colon B^{*} \times B \to \La$ denotes the evaluation pairing.
\end{proof}
 
 \begin{lemma} \label{lem:reduced hermitian}
 The form $(L_M, h_M)$ is hermitian.
\end{lemma}
\begin{proof}  Since $s_M$ is hermitian,  we have the formula $\ad s_M = (\ad s_M)^* \circ \gamma$ from Lemma \ref{lem:hermitian}. 
After dualizing we obtain $(\ad s_M)^* = \gamma^* \circ (\ad s_M)^{**} $. Note that $\gamma^* \colon B^{***} \to B^*$ is the left  inverse of $\gamma_1\colon B^* \to (B^{**})^{*}= (B^{*})^{**}$,  given by $v \mapsto v^{**}$, for any right $\La$-module $B$. This follows from the calculation
$$\la (\gamma^*\circ \gamma_1)(v), b\ra =\la \gamma^*(v^{**}), b\ra= \la v^{**}, \gamma(b)\ra = v^{**}(b^{**}) = \overline{b^{**}(v)} = v(b),$$
for all $b \in B$ and all $v \in B^*$.  Therefore
 $\gamma_1\circ (\ad s_M)^* = (\ad s_M)^{**}$. 

By definition, $\ad h_M\circ (\ad s_M)^* = \id$, so we have $(\ad s_M)^{**} \circ (\ad h_M)^* = \id$ by dualizing. Therefore, $(\ad h_M)^* \circ \gamma_1 = \ad h_M$.  
From Lemma \ref{lem:hermitian} applied to $B = \pi_2(M)^*$ and $h = h_M$, we conclude that $h_M$ is hermitian symmetric.
\end{proof}
\begin{remark}
If $\pi = \pi_1(M)$ is any finitely presented group with tame cohomology, then the splicing argument above shows directly that the  dual of the equivariant intersection form $s_M$ induces a well-defined, non-singular hermitian form $h_M \colon L_M \times L_M \to \La$  on $\pi_2(M)^*$. However, if $\cd \pi \geq 4$, the module $\pi_2(M)^*$ might not be stably free (take $\pi = \bZ^4$).  
\end{remark}

\begin{definition} [\spinp manifold]\label{def:reduced SW}
 If $\cd \pi_1(M)\leq 3$, we define the \emph{reduced Stiefel-Whitney class} 
$$w_M\colon L_M \to \cy 2$$ by the formula
$w_M(x) = \varepsilon_1 (h_M(x,x)) \pmod 2$, for all $x\in L_M$. Here $\varepsilon_1\colon \La \to \bZ$ is the coefficient of $\lambda \in \La$ at the identity element. The form $h_M$ is even if and only if $w_M \equiv 0$. In this case, we say that $h_M$ has \emph{even} reduced $w_2$-type, and that $M$ is a \emph{\spinp} manifold.
\end{definition}
We are indebted to Peter Teichner for pointing out that $h_M$ is not always an even form. 
\begin{remark}  Examples of spin $4$-manifolds with even reduced $w_2$-type  were given by Plotnick \cite[\S 3]{Plotnick:1986} using a ``twist-spinning" construction applied to any closed, oriented, aspherical $3$-manifold (compare \cite[\S 6.2]{Kasprowski:2017}). In these examples, $s_M$ is the hyperbolic form on the augmentation ideal $I(\pi) \subset \Zpi$. \end{remark}

We conclude this section by showing that the thickened double construction provides a \spinp manifold for every fundamental group $\pi$, whenever $\cd \pi \leq 3$.

Note first that the module $\pi_2(M)^{**}$ supports the non-singular hermitian form 
$\widehat s_M (u,v) = \la(\ad s_M)^{**}(u), v\ra$ whose adjoint is $(\ad s_M)^{**}$. The hermitian property follows as  in Lemma \ref{lem:hermitian}.

\begin{lemma}
If $s_M$ is strongly even, then $h_M$ has even reduced $w_2$-type.  
\end{lemma}
\begin{proof}
We first show that the homomorphism $\widehat w_M \colon \pi_2(M)^{**} \to \cy 2$  associated to  the form $\widehat s_M$ by the formula
$u \mapsto \varepsilon_1(\widehat s_M(u,u)) \pmod 2$ is determined by $w_M$. To check this we will write $u = \ad h_M(x)$, for  a unique $x \in \pi_2(M)^*$, and use the formula
$\gamma_1\circ (\ad s_M)^* = (\ad s_M)^{**}$ from the proof of Lemma \ref{lem:reduced hermitian}. Now we have:
\begin{align*}
 \la (\ad s_M)^{**} (u), u\ra &= \la (\ad s_M)^{**} (\ad h_M(x)), \ad h_M(x)\ra\\
 &= \la \gamma_1\circ (\ad s_M)^* (\ad h_M(x)), \ad h_M(x))\ra \\&= \la \gamma_1(x) , \ad h_M(x)\ra = \ad h_M(x)(x),
\end{align*}
for all $u \in \pi_2(M)^{**}$, with $u = \ad h_M(x)$.
Therefore  $\widehat w_M (u) = w_M(x)$,  as required. 

Now suppose that $s_M = \lambda + T\lambda$ for some even hermitian form
 $(\pi_2(M), \lambda)$. We have the formula
 $$(\ad s_M)^{**} = (\ad\lambda)^{**} + (\ad \lambda)^{***} \circ \gamma^{**}$$
 and we claim that 
 $$((\ad \lambda)^{***} \circ \gamma^{**})(u)(v) = \overline{((\ad\lambda)^{**} )(v)(u)}.$$
 This follows since 
 \begin{align*}
 \la (\ad \lambda)^{***} \circ \gamma^{**})(u), v\ra &= \la (\ad \lambda)^{***} (\gamma^{**}(u)), v\ra = \la \gamma^{**}(u), (\ad \lambda)^{**}(v)\ra\\
 &=\la u, \gamma^*((\ad \lambda)^{**}(v))\ra =\la (\ad \lambda)^{**}(v), \gamma(u)\ra \\
 &= \overline{\la (\ad\lambda)^{**}(v), u\ra}
 \end{align*}
 for all $u, v \in \pi_2(M)^{**}$. Therefore $s_M$ strongly even implies $\widehat w_M = 0$,  and hence $w_M =0$.
\end{proof}

\begin{corollary}\label{cor:allgroups}
 Let $\pi$ be a finitely presented group with $\cd\pi \leq 3$. Then there exists a smooth, closed, oriented \spinp manifold $M^4$ with $\pi_1(M) = \pi$.
\end{corollary}
\begin{proof} Let $M$ be a smooth, spin,  thickened double $4$-manifold as constructed in Lemma \ref{lem:double}. Then the equivariant intersection form $s_M$ is strongly even. This follows from the formula in \eqref{eq:doubleform}, since the metabolic form $s_M$ is determined by an even hermitian form $(\pi_2(K), \delta_K)$ which depends on the the particular thickening. However, 
$\pi_2(K)$ may be assumed to be a free $\La$-module by stabilization, and therefore $\delta_K$ is strongly even.
\end{proof}

\section{The reduced $2$-type $P$}\label{sec:seven}

In order to prove Theorem A, we will follow the strategy of \cite[\S  2]{HKT09} involving the reduced normal $2$-type and the ``modified" surgery theory developed by Kreck \cite{Kreck99}. Since we have restricted our attention to  spin topological  manifolds, the \emph{reduced normal $2$-type} for a given $4$-manifold $M$ is a fibration $B \to BSTOP$, with total space
$$B := B(M) = P \times BTOPSPIN$$
 defined by the second factor projection  $P \times BTOPSPIN \to BTOPSPIN$ and the natural map $BTOPSPIN \to BSTOP$. Here we are using the stabilized classifying spaces for topological $\bbR^n$-bundles, with $w_1=0$ for $BSTOP$ and $w_1, w_2 = 0$ for $BTOPSPIN$ (see \cite[p.~190]{kirby-siebenmann1}). 

 It remains to describe the reduced $2$-type $P = P(M)$.

\begin{definition}\label{def:reduced 2type}
 Let $M$ be a closed, oriented, spin 4-manifold with fundamental group $\pi$. We define the \emph{reduced $2$-type of $M$} as follows:  let $P= P(M)$ be the two-stage Postnikov system with $\pi_1(P) = \pi_1(M)=\pi$ and $\pi_2(P)=\LM=\pi_2(M)^*$.  The total space $P$ is defined by a fibration over $K(\pi,1)$ with fibre $K(\LM,2)$ and classified by $k_P$ in $H^3(\pi;\LM)$:
\begin{equation*}
\xymatrix@R-10pt@C-5pt{ K(\LM,2) \ar[r] & P \ar[d] & \\ & K(\pi,1) \ar[r] & H^3(\pi; \LM) }
\end{equation*}

We have a reference map $c_M\colon M\to P$ which factors through the algebraic 2-type of $M$ (classified by $\pi_1(M)=\pi$, $\pi_2(M)$, and $k_M \in H^3(\pi;\pi_2(M))$). The $k$-invariant $k_P$ is given by the image of $k_M$ under the  map induced by  $\pi_2(M) \to \pi_2(P)$.
\end{definition}

The map  $\pi_2(M) \to \pi_2(P)$ is neither an injection nor a surjection in general because $P$ is given by the union of cells of dimension $\geq 2$ with $M$. In particular, the map $c_M$ is $1$-connected but not $2$-connected.

\begin{definition} For $M$ a closed, spin, TOP $4$-manifold, 
a map $c\colon M\to P$ to a space $P$ is called a \emph{reduced 3-equivalence} if $c$ induces an isomorphism $c_*\colon \pi_1(M) \cong \pi_1(P)$ and an isomorphism $c^*\colon \pi_2(P)^* \cong \pi_2(M)^*$. If $c$ is a reduced 3-equivalence, then $c \times \nu_M\colon M \to P \times BTOPSPIN$ is called a \emph{reduced normal 2-smoothing}.
\end{definition}

We have a result similar to \cite[Proposition 2.11]{HKT09}.
\begin{proposition}\label{prop:seventhree}
 Let $M$ and $M'$ be closed, oriented, spin, TOP $4$-manifolds with fundamental group $ \pi$. Suppose that there is an isometry $\fQ(M) \cong \fQ(M')$ given by isomorphisms $\alpha\colon \pi_1(M) \to \pi_1(M')$ and $\beta\colon \pi_2(M) \to \pi_2(M')$. Then there is a $3$-coconnected fibration $P\times BTOPSPIN \to BSTOP$ admitting  two reduced normal 2-smoothings $M \to P \times BTOPSPIN$ and $M' \to P \times BTOPSPIN$ that induce $(\alpha, \beta)$. 
\end{proposition}

\begin{proof} After identifying $\pi_1(M)=\pi$ and $\pi_1(M')$ via the given isomorphism $\alpha$, the map $\beta\colon \pi_2(M) \to \pi_2(M')$ is $\pi$-equivariant, and 
$(\ad s_M) = \beta^* \circ (\ad s_{M'})\circ \beta $. If we let $P = P(M)$ and $P'= P(M')$ denote the reduced 2-types for $M$ and $M'$ respectively, then we have a fibre homotopy equivalence
$P \to P'$ over $K(\pi, 1)$ induced by the given isometry $(\alpha, \beta)$  preserving the $k$-invariants.
\end{proof}

\section{The cohomology of $P$}\label{sec:eight}

In order to prove Theorem A, we need to calculate the topological spin bordism group $\Omega_4^{Spin}(P)$,  and this requires some information about $H^*(P;\La)$ and $H_*(P;\bbZ)$. 

If $c\colon M \to P$ is a reduced $3$-equivalence,  we have the same 4-term exact sequence for $P = P(M)$ as we do for $K$ and $M$. The following diagram commutes: 
\eqncount
\begin{equation}\label{4term for P} 
\vcenter{\xymatrix@C-5pt{0 \ar[r]& H^2(\pi; \La) \ar[r]\ar@{=}[d]& H^2(P; \La) \ar[r]^-{eval}\ar[d]_{\cong}^{c^*}&
\Hom_\La(\pi_2(P),\La) \ar[r]^-{d_3}\ar[d]_{\cong}^{c^*}&H^3(\pi;\La) \ar[r]\ar@{=}[d]& 0\cr
0 \ar[r]& H^2(\pi; \La) \ar[r]& H^2(M; \La) \ar[r]^-{eval}&
\Hom_\La(\pi_2(M),\La) \ar[r]&H^3(\pi;\La) \ar[r]& 0}}
\end{equation}
and we see that $H^2(P;\La)\cong H^2(M;\La)$.
Here is a partial calculation of the cohomology of the reduced $2$-type $P = P(M)$.

\begin{figure}
\begin{center}
\begin{tikzpicture}
[x=32pt, y=32pt, line width=1pt]
\draw[black!60, step=32pt, line width=.3pt] (-.15,-.15) grid (3.5, 3.5);
\draw (-.15, 0) -- (3.5, 0) node[anchor=west] {$p$};
\draw[-] (0, -.15) -- (0, 3.5) node[anchor=south] {$q$};
 \path (3,0) node[anchor=north] {$H^3(\pi;\La)$};
\path (0,2) node[anchor=east] (q) {$\Hom_\La(H_2(P;\La),\La)$};
\draw[-stealth, red, line width=2pt, domain = 0:2.9]
	plot (\x,{-2/3*\x+2});
\foreach \i in {0,...,3}{
	\draw[black, fill=white, line width=1pt] (\i,1) circle (3pt) (\i,2) circle (3pt) (\i,3) circle (3pt);}
\foreach \i in {1,2,3}{
	\fill[black]
		(\i,0) circle (1pc/3.5);}
\fill[black]
	(0,2) circle (1pc/3.5);
\end{tikzpicture} 
\caption{The $E_2$-page of the spectral sequence $E_2^{p,q} = \Ext_{\La}^p(H_q(P;\La),\La) $  converging to $H^*(P;\La)$. The solid dots represent possible  non-zero terms.}
\label{E2 page cohomology of P}
\end{center}
\end{figure}
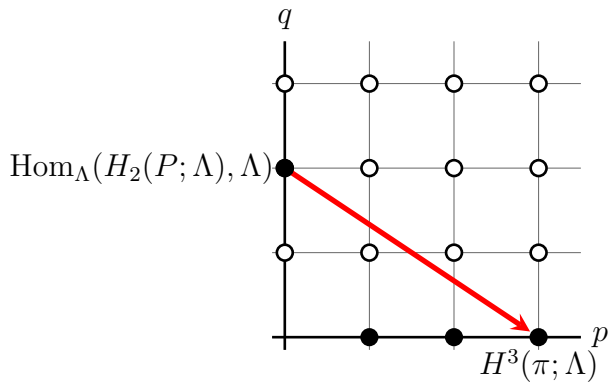

\begin{lemma} \label{cohomology of P} 
 $H^1(P;\La)\cong H^1(\pi;\La)$, 
 $H^2(P;\La)$ is an extension of $\ker d_3\subseteq \pi_2(P)^*$ by $H^2(\pi;\La)$, and  $H^3(P;\La)=0$.\end{lemma}
\begin{proof}
We  look at the universal coefficient spectral sequence $E^{p,q}_2 \cong \Ext^p_\La(H_q(P;\La),\La)$ to compute $H^*(P;\La)$. For $q=0$ we have 
$E^{p,0}_2 \cong \Ext^p_\La(\bZ,\La) = H^p(\pi;\La)$. When $q$ is odd, $H_q(P;\La) \cong H_q(\wP;\bbZ) =0$ since $\wP$ is a product of copies of $\CP^\infty$.  Thus the $E^{p,1}_2$- and $E^{p,3}_2$-terms are zero. Additionally, since $H_2(P;\La) \cong \pi_2(P)$ is stably $\La$-free (by Corollary \ref{cor:stableiso} and Lemma \ref{s_M star}), the terms $\Ext^p_\La(H_2(P;\La),\La) = 0$ for $p>0$.
The $E^{1,2}_2$-term of the spectral sequence will therefore be zero. 
The only non-zero differential is the $d_3$ map in Figure \ref{E2 page cohomology of P}, which is surjective by (\ref{4term for P}). Hence the $E^{3,0}_\infty$-term is zero, and so $H^3(P;\La)=0$. 
\end{proof}

\begin{lemma}\label{no obstruction}
For a closed, oriented, TOP 4-manifold $M$ with reduced 2-type $P$, $$H^3(P,M;\La) = 0.$$ 
\end{lemma}
\begin{proof}
We use the long exact sequence in cohomology of the pair $(P,M)$ with $\La$-coefficients. We have the isomorphisms $H^3(M;\La) \cong H_1(M;\La) \cong H_1(\wM;\bbZ) = 0$, so the long exact sequence becomes 
\begin{equation*}
\dotsm \to H^2(P;\La) \to H^2(M;\La) \to H^3(P,M;\La) \to H^3(P;\La) \to 0.
\end{equation*}

From (\ref{4term for P}), the isomorphism $H^2(P;\La)\cong H^2(M;\La)$ forces the middle map in the above sequence to be zero. This yields the isomorphism $H^3(P,M;\La)\cong H^3(P;\La)$, which is zero by Lemma \ref{cohomology of P}.
\end{proof}

\section{The homology of $P$}\label{sec:nine}

In this section, let $M$ be a closed, spin, TOP $4$-manifold such that $\pi_1(M) = \pi$ is a finitely presented group with $\cd \pi\leq 3$. 
To compute the homology of the reduced 2-type $P=P(M)$, we use the Serre spectral sequence for the fibration $\wP \to P \to K(\pi,1)$. This spectral sequence with integral coefficients has the $E^2$-page
$$E^2_{p,q} = H_p(\pi; H_q(\wP))$$
and we only need the homology of $P$ up to dimension 5. 
Note that $\wP$ is a restricted product of copies of $\CP^\infty$, so $H_q(\wP) = 0$ for $q$ odd. We have already seen that $H_2(\wP) = \pi_2(P)$ is a stably free $\La$-module. Therefore, by \cite[Lemma 2.2]{HambletonKreck88}, we see that $H_4(\wP) = \Gamma(\pi_2(P))$ is also stably free\footnote{The proof given for this lemma  applies without change to infinite groups.}. The notation $\Gamma(\pi_2(P))$ refers to Whitehead's quadratic functor, defined in \cite[\S 5]{HambletonKreck88}.
In summary:
\begin{enumerate}
\item $E^2_{p,0} \cong H_p(\pi; \bbZ)$, which is zero when $p\geq 4$,
\item  $E^2_{0,q} \cong H_0(\pi;H_q (\wP))\cong H_q(\wP) \otimes_\La \bbZ$,
\item  $E^2_{p,q}=H_p(\pi; H_q(\wP))=0$ is zero for odd $q$,
\item $E^2_{p,2}=H_p(\pi;H_2(\wP))=0$, for $p>0$,
\item $E^2_{p,4}=H_p(\pi; H_4(\wP) )= H_p(\pi;\Gamma(\pi_2(P))) = 0$, for $p>0$.
\end{enumerate}

In the $E^2$-page, all $d^2$ maps that affect $H_i(P)$, $i\leq 5$, are zero. In the spectral sequence, the only possibly non-zero  differential is 
$d^3\colon H_3(\pi) \to H_0(\pi; \pi_2(P))$.
\begin{proposition}\label{d3 P}
$d^3\colon  H_3(\pi) \to H_0(\pi; \pi_2(P))$ is injective. 
\end{proposition}
The first step of the proof is to establish this result for the  minimal model $M_0$.

\begin{lemma}\label{d3 P minimal} Let $P = P(M_0)$ be the reduced $2$-type of the minimal model. Then the map
$d^3\colon  H_3(\pi) \to H_0(\pi; \pi_2(P))$ is injective. 
\end{lemma}

\begin{proof}
The injectivity argument comes from comparing the same $d^3$ maps in three  spectral sequences, the first of which is for $H_*(K)$. In the spectral sequence converging to $H_*(K)$, since $H_3(K)$ surjects onto the $E^\infty_{3,0}$ term and $H_3(K)=0$, the differential $$d^3\colon H_3(\pi)\to H_0(\pi;\pi_2(K))$$ must be injective. 
Recall that $M_0$ is the double of a thickening $Y:=Y(K)$ of $K$, and we view $M_0$ as the boundary of $Y\times I$; this gives a map $M_0 \hookrightarrow Y\times I \simeq K$. The reduced 2-type $P$ is constructed by attaching cells of dimension 2 and higher to $M_0$. We define 
$$P^{(3)} := M_0 \cup \bigcup_{\alpha}e^2_\alpha \cup \bigcup_{\beta} e^3_\beta,$$
 as the union of $M_0$ with only the  2-cells and 3-cells from $P$. Since $H^3(P,M_0;\pi_2(K))=0$ by Lemma \ref{no obstruction}, obstruction theory tells us the map $M_0 \to K$ extends over $P^{(3)}$, and we obtain an induced map $H_0(\pi;\pi_2(P^{(3)})) \to H_0(\pi; \pi_2(K))$. By commutativity of the diagram below, the $d^3$ map in the spectral sequence converging to $H_*(P^{(3)})$ must also be injective.
\begin{equation*}
\xymatrix{H_0(\pi;\pi_2(K)) & \ar[l] H_0(\pi;\pi_2(P^{(3)})) \\ H_3(\pi) \ar@{=}[r] \ar@{^(->}[u]^{d^3} & H_3(\pi) \ar[u]^{d^3}}
\end{equation*}

It remains to compare the $d^3$ differentials for $H_*(P)$ and $H_*(P^{(3)})$. 
We claim that $\pi_2(P) \cong \pi_2(P^{(3)})$: the relative homologies $H^i(P,P^{(3)};\La)$ vanish in dimension $i=2,3$, so the isomorphism is given by the long exact sequence of the pair. The injectivity of $d_3\colon  H_3(\pi) \to H_0(\pi;\pi_2(P^{(3)}))$ implies that $d^3\colon H_3(\pi) \to H_0(\pi;\pi_2(P))$ is also injective. 
\end{proof}

The following lemma  is used  in the proof of Proposition \ref{d3 P}.
\begin{lemma} \label{stabilizing with X} Let $M$ be a closed, oriented, TOP $4$-manifold with $\pi_1(M) = \pi$, and let $X$ be a closed, simply connected $4$-manifold. Then the map $d^3\colon  H_3(\pi) \to H_0(\pi; \pi_2(P(M)))$ is injective if and only if the map $d^3\colon  H_3(\pi) \to H_0(\pi; \pi_2(P(M\# X)))$ is injective.
\end{lemma}
\begin{proof}
We begin by comparing $M$ and $M\#X$. By removing the top dimensional cells of $M$ and $M\# X$, we get an inclusion $M^o \hookrightarrow (M\# X)^o$, the latter of which is just $M^o$ wedged with a collection of $n$ 2-spheres arising from $X^o$. This inclusion induces a split injection $\pi_2(M^o) \to \pi_2((M\# X)^o)$, and so $\pi_2(M)$ is stably isomorphic to $\pi_2(M\#X)$:
$$ \pi_2(M\#X) \cong \pi_2((M\#X)^o) \cong \pi_2(M^o) \oplus \La^n \cong \pi_2(M) \oplus \La^n,$$
where $n = b_2(X)$.
If $P(M)$ is the reduced 2-type of $M$, and $P(M\# X)$ is the reduced 2-type of $M\#X$, then it follows that
$\pi_2(P(M))$ is stably isomorphic to $ \pi_2(P(M\#X))$, and therefore
$H_0(\pi;\pi_2(P(M))) \cong H_0(\pi;\pi_2(P(M\# X)))$. The conclusion about injectivity for the maps $d^3$ now follows by naturality of the spectral sequences with respect to the map $P(M) \to P(M \# X)$.
\end{proof}

\begin{proof}[The proof of Proposition \ref{d3 P}]
By Proposition \ref{prop:stable},  we have
$$ M\#\CP^2\#\overline{\CP^2} \# r(S^2 \times S^2) \approx M' \#s(S^2 \times S^2).$$
where $M' = M_0 \#p \CP^2 \#q \overline{\CP^2}$ is a suitable stabilization of the minimal model $M_0$. 

This homeomorphism allows us to equate their second homotopy groups, and thus identify their reduced 2-types. By applying Lemma \ref{stabilizing with X} several times, the $d^3$ map for $P(M)$ is injective if and only if 
the $d^3$ map for $P(M_0)$ is injective, and Lemma \ref{d3 P minimal} completes the proof.
\end{proof}

\begin{remark}
The same $d^3$ map in the spectral sequence converging to $H_*(M)$ is injective as well, given by naturality of the spectral sequences under the map $M\to K$. \end{remark}

\begin{proposition} The integral homology  $H_i(P)$, for $i \leq 5$,  is given as follows:
\begin{enumerate}
\item $H_0(P) = \bZ$ and $H_1(P) = H_1(\pi)$. 
\item $H_2(P)$ is an extension of  $H_2(\pi)$ by $\coker \{d^3\colon H_3(\pi)  \to H_2(\widetilde P) \otimes_{\La} \bZ$\}.
\item $H_3(P) =   H_5(P) = 0$
\item $H_4(P) \cong H_4(\wP) \otimes_\La \bbZ$.
\end{enumerate}
In addition, $H_3(P;\cy 2) =0$.
\end{proposition}
\begin{proof}
The table summarizes the calculations above. For \RAAGs\ the 
 extension in (ii) is split. The argument showing that $H_3(P;\cy 2) =0$ follows that same steps as above. The details will be left to the reader. 
\end{proof}

\section{The spin bordism groups $\Omega^{Spin}_*(P)$  }
\label{sec:ten}

The Atiyah-Hirzebruch  spectral sequence is used to compute the topological (or smooth) spin bordism groups of the reduced 2-type $P=P(M)$, denoted as above by $\Omega^{Spin}_*(P)$. 
The $E^2$-page of the spectral sequence is given by $H_p(P;\Omega^{Spin}_q(\pt))$, and the relevant Spin bordism groups of a point are given below:
\begin{equation*}
\Omega^{Spin}_q(\pt) = \bbZ, \cy 2, \cy 2, 0, \bbZ\quad  \textrm{for } q = 0,1,2,3,4
\end{equation*}
The topological case differs from the smooth case only in the divisibility of the signature, so that $\Omega^{TopSpin}_4(pt)  = 8\bZ$ while $\Omega^{Spin}_4(pt) = 16\bZ$.  
\begin{figure}[h]
\begin{center}
\begin{tikzpicture}
[x=32pt, y=32pt, line width=1pt]
\draw[black!60, step=32pt, line width=.3pt] (-.05,-.25) grid (5.5, 4.5);
\draw (-.05, 0) -- (5.5, 0) node[anchor=west]{$p$};
\draw[-] (0, -.25) -- (0, 4.5) node[anchor=south]{$q$};
\clip (-1,-.25) rectangle (5.5, 4.5);
\draw[stealth-, red, line width=2pt, domain = .1:2]
	plot (\x,{-1/2*\x+1});
\draw[stealth-, red, line width=2pt, domain = .1:2]
	plot (\x,{-1/2*\x+2});
\draw[stealth-, red, line width=2pt, domain = 2.1:4]
	plot (\x,{-1/2*\x+2});
\draw[stealth-, red, line width=2pt, domain = 2.1:4]
	plot (\x,{-1/2*\x+3});
\foreach \j in {0,...,4}{
	\draw[black, fill=white, line width=1pt] (3,\j) circle (3pt) (5,\j) circle (3pt);}
\foreach \i in {0,...,5}{
	\draw[black, fill=white, line width=1pt] (\i,3) circle (3pt);}
\draw[black, fill=white, line width=1pt] (5,0) circle (3pt);
\foreach \i in {0,4}{\path (0,\i) node[anchor=east] (z\i) {$\bbZ$};}
\foreach \i in {1,2}{\path (0,\i) node[anchor=east] (z2\i) {$\bbZ/2$};}
\path (0,3) node[anchor=east] (0) {$0 \ $};
\foreach \i in {0,1,2}{
	\foreach \j in {0,1,2}{
		\fill[black]
			(\i,\j) circle (1pc/3.5);
		}}
\fill[black]
	(4,0) circle (1pc/3.5)
	(4,1) circle (1pc/3.5)
	(0,4) circle (1pc/3.5);
\end{tikzpicture}
\caption{The $E^2$-page of the spectral sequence converging to $\Omega^{Spin}_*(P)$.} 
\label{bordism sseq}
\end{center}
\end{figure}
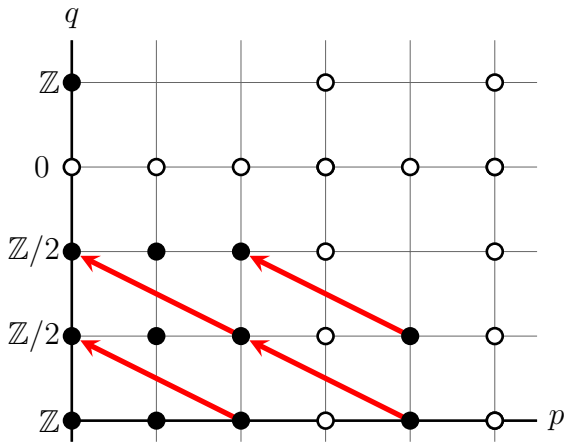

In Figure \ref{bordism sseq}, we have included the information about $H_*(P)$ from the last section. The two $d^2$ maps on the left in the $E^2$-page  are both zero,  otherwise $\Omega^{Spin}_*(\pt)$ would not split off in the $E^\infty$-page of $\Omega^{Spin}_*(P)$. The other two $d^2$ maps are the duals of the $Sq^2$ maps composed with reduction mod 2. 

Consider the commutative diagram below that arises from the fibration $\wP \to P \to K(\pi,1)$. We take homology with $\cy 2$-coefficients.
\begin{equation*}
\xymatrix{H_4(\wP) \ar@{->>}[r] \ar@{->>}[d] & H_2(\wP) \ar[d] & & \\
H_4(P) \ar[r]^{DSq^2} \ar[d] & H_2(P) \ar[r] \ar@{->>}[d] & \cok(DSq^2) \ar[r] \ar@{.>}[d] & 0 \\
0=H_4(\pi) \ar[r] & H_2(\pi) \ar@{=}[r] \ar[d] & H_2(\pi) & \\
& 0}
\end{equation*}
The above map labeled $DSq^2$ is the dual of  $Sq^2\colon H^2(P) \to H^4(P)$. Since $\wP$ is a restricted product of copies of $\CP^\infty$ and inverse limits are left exact, $H^2(\wP)$ injects into $H^4(\wP)$, and thus $H_4(\wP)$ surjects onto $H_2(\wP)$. But $H_4(\pi)=0$ by assumption,  so 
by exactness of the middle vertical sequence, $\cok(DSq^2) \cong H_2(\pi; \cy 2)$. 
Since $H_5(P;\bZ) = 0$, the term $H_2(\pi; \cy 2)$ survives to $\Omega_4^{Spin}(P)$.

\begin{proposition}\label{prop:bordism}
 The topological or smooth   spin bordism groups of the reduced normal $2$-type $P=P(M)$ for $M$ are detected by an injection
$$\Omega^{Spin}_4(P) \subseteq \bZ  \oplus H_2(\pi;\cy 2) \oplus H_4(P;\bZ).$$
The invariants are the signature,  an invariant in $H_2(\pi; \cy 2)$, and the fundamental class $c_*[M] \in H_4(P; \bZ)$.
\end{proposition}
\begin{proof}
This follows from our discussion of the differentials in the spectral sequence.
\end{proof}
The elements of  $\Omega^{Spin}_4(P)$ are represented by pairs $(N, f)$, where $N$ is a closed, spin $4$-manifold and $f \colon N \to P$ is a reference map. We will always assume that $f$ induces an isomorphism of fundamental groups.

For \spinp manifolds (recall Definition \ref{def:reduced SW}), we will show that the bordism invariant in $H_2(\pi;\cy 2)$ is determined by the other invariants. The method follows  \cite[\S 5]{HKT09}, where the authors define a subset $\Omega_4(P)_M \subset \Omega^{Spin}_4(P)$, called the \emph{normal structures}. We need to adapt this notion to our setting.

\begin{definition} Let $M$ be a \spinp manifold.
The set of \emph{normal structures} $\Omega_4(P)_M \subset \Omega_4^{Spin}(P)$ consists of the spin bordism classes $(N,f)$, where $N$ has even reduced $w_2$-type,  such that  $f_*[N] = c_*[M]$ and $\sign(N) = \sign(M)$.
\end{definition}
We note this set is non-empty, since $(M,c) \in \Omega_4(P)_M$ whenever $M$ is a \spinp manifold. In Lemma \ref{lem:tenthree} we will check that the definition of a normal structure is independent of the choice of representative $(N, f)$ up to spin bordism.

\medskip
Here is a useful observation. 

\begin{lemma}\label{lem:tenfour}  Let $(N,f)$ represent an element of $\Omega^{Spin}_4(P)$. 
If $f_*[N] = c_*[M]$, then the reduced intersection form $(\pi_2(N)^*, h_N)$ restricted to the image of $f^*\colon \pi_2(P)^* \to \pi_2(N)^*$ is isometric to  $(\pi_2(M)^*, h_M) $.
\end{lemma}

\begin{proof}
Since $f_*[N] = c_*[M] \in H_4(P;\bZ)$, we have the following commutative diagram:
\eqncount
\begin{equation}\label{eq:pullback}
\vcenter{\xymatrix@R+5pt@C+15pt{
\pi_2(M) \ar[d]^{c_*} & H^2(M;\La) \ar[l]^{\cap [M]}_{\approx} \ar[r]^{eval}& \pi_2(M)^*\\
\pi_2(P) & H^2(P;\La) \ar[l]^{\cap c_*[M]} \ar[r]^{eval}\ar[d]^{f^*}\ar[u]_{c^*}^{\approx}
& \pi_2(P)^*\ar[d]^{f^*}\ar[u]^{c^*}_{\approx}\\
\pi_2(N) \ar[u]^{f_*} & H^2(N;\La) \ar[l]^{\cap [N]}_{\approx} \ar[r]^{eval}& \pi_2(N)^*\\
}}
\end{equation}
The composite on the top row defines $\ad s_M$, and the composite on the
bottom row defines $\ad s_N$. After dualizing each term  in the diagram, each square still commutes,  and the upper vertical maps ($c_*$ and $c^*$) as well as all the horizontal maps become isomorphisms.
 The composites in the top and bottom rows of the dualized diagram give $(\ad s_M)^*$ and
$(\ad s_N)^*$,  respectively.  

Therefore the inverses of the dualized top and bottom row composites give $\ad h_M$ and $\ad h_N$, respectively.  One can check by a diagram chase that the inverse of the dualized composite in the  middle row is isometric via $(c^*, c^{**})$ to the adjoint of
$h_M$. A similar diagram chase  shows that the inverse of this middle row composite is also isometric to the pull-back of the adjoint of $h_N$ via $(f^*, f^{**})$.
\end{proof}

\begin{corollary}\label{cor:reducediso}
If $f\colon N \to P$ is a reduced $3$-equivalence and $f_*[N] = c_*[M]$, then 
$(\pi_2(M)^*, h_M) \cong (\pi_2(N)^*, h_N) $. 
\end{corollary}
\begin{proof} We have the isomorphisms:
$$(\pi_2(N)^*,h_N) \cong (\pi_2(P)^*,(h_M)^*) \cong  (\pi_2(M)^*,h_M)$$
induced by the maps in the diagram \eqref{eq:pullback}.
\end{proof}

We will now define a map
 $\theta_M\colon \Omega_4(P)_M \to L_4(\Zpi)$ as in \cite[Definition 5.8]{HKT09}. 
Let $(N,f)$ represent an element of $\Omega^{Spin}_4(P)$.   
After preliminary surgeries we may assume that the map $f$ is 2-connected, so that 
 $(\pi_2(N)^*,h_N)$ is stabilized by direct sum with a hyperbolic form $H(\La^r)$. Note that if $N$ is a \spinp manifold, then the stabilized manifold also has even reduced $w_2$-type. Let
$$V:= \ker \left ( f_*\colon \pi_2(N) \to \pi_2(P)\right ).$$
Since $\pi_2(P)$ is a stably free $\Lambda$-module and $f_*$ is surjective, this sequence splits.
 After applying $\Hom_\La( -, \La)$, we obtain a (split) short exact sequence
\eqncount
\begin{equation}\label{eq:firstdual}
 0 \to \pi_2(P)^* \xrightarrow{f^*} \pi_2(N)^* \to V^*\to 0.
 \end{equation}
Since $\pi_2(N)^*$ is stably free, so is the module $V^*$. 

\begin{definition} Suppose that $(N,f)$ has $2$-connected reference map. 
Let $(V^*,\lambda_{N,f})$ denote the restriction of the form $(\pi_2(N)^*,h_N)$ to the orthogonal complement of the submodule $\pi_2(P)^*$. 
\end{definition}
\begin{lemma}\label{lem:tenseven}
If $M$ is a \spinp manifold and $(N,f) \in \Omega_4(P)_M$ is a normal structure with $2$-connected reference map,  then $(V^*,\lambda_{N,f})$ is a non-singular even  hermitian form on a finitely generated, stably free $\Lambda$-module.
\end{lemma}
\begin{proof} 
By Lemma \ref{lem:tenfour}, the form $h_N$ restricted to the image of $f^*$ is isometric to $h_M$. Since $h_M$ is non-singular, we obtain an isometric splitting
$$(\pi_2(N)^*,h_N) \cong (\pi_2(M)^*,h_M)\oplus(V^*,\lambda_{N,f})$$
and therefore $(V^*,\lambda_{N,f})$ is a non-singular form. Since $h_N$ is even by assumption,  the form $(V^*,\lambda_{N,f})$ is also even.
\end{proof}

 We now return to detecting the bordism invariant in $H_2(\pi;\cy 2)$ for the elements  $(N,f) \in \Omega_4(P)_M$.  Let 
$(N,f)$ be a normal structure representing an element of  $ \Omega_4(P)_M$, where $M$ is a \spinp manifold.  Lemma \ref{lem:tenseven} provides 
 a non-singular \emph{even} hermitian form  $(V^{*}, \lambda_{N,f})$ on a finitely generated, stably free $\La$-module. In particular, the form $\lambda_{N,f}$ admits a unique quadratic refinement. We define
$\theta_M\colon \Omega_4(P)_M \to \tilde L_4(\Zpi)$
by  setting
$$ \theta_M(N,f) = (V^{*}, \lambda_{N,f}) \in \tilde L_4(\Zpi), $$
where $\tilde L_4(\Zpi)$ denotes the reduced surgery obstruction group (represented by quadratic forms with signature zero), for surgery up to homotopy equivalence. 

\begin{lemma} \label{lem:tenthree} Let $M$ be a \spinp manifold.
 The map $\theta_M$ is well-defined. In addition, $\theta_M(N, f)  = 0$ if $(N,f)$ is spin bordant over $P$ to a reduced $3$-equivalence.
\end{lemma}
\begin{proof}   If $(N, f)$ and $(N', f')$ are spin bordant over $P$ with $2$-connected reference maps, then $N$ and $N'$ are stably homeomorphic over $P$ by \cite[Corollary 3]{Kreck99}. However,  an arbitrary element $(N,f) \in \Omega_4^{Spin}(P)$ is spin bordant to one with $2$-connected reference map,  with the property that $s_N$ is stabilized by a hyperbolic form. Hence the parity (even or odd) of $h_N$ is  preserved by this operation.  

It follows that any element $(N', f')$ which is spin bordant to $(N,f) \in \Omega_4(P)_M$ is also a normal structure. 
The stable homeomorphism argument now shows that $\theta_M$ is well-defined (compare with\cite[Lemma 5.9]{HKT09}). 

Now suppose that $f\colon N\to P$ is a reduced 3-equivalence, and that $(N', f') \in \Omega_4(P)_M$ is a $2$-connected representative in the same spin bordism class.  We have 
$$f_*[N] = f'_*[N'] = c_*[M] \in H_4(P).$$
 Lemma \ref{lem:tenfour} provides the isometric direct sum splittings: 
$$ (\pi_2(N)^*,h_N)\oplus H(\La^r) = (\pi_2(N')^*,h_{N'}) \cong (\pi_2(M)^*,h_M)\oplus(V^*,\lambda_{N,f})$$
and the isomorphisms:
$$(\pi_2(N)^*,h_N) \cong (\pi_2(P)^*,(h_M)^*) \cong  (\pi_2(M)^*,h_M)$$
 which give a relation in the $L$-group, since all these forms are non-singular and even hermitian. It follows that $[(V^*, \lambda_{N,f})] = [H(\La^r) ]= 0 \in L_4(\Zpi)$, and $\theta_M(N,f) = 0$. 
\end{proof}

The next step is to define a map 
$$\rho_M \colon \Omega_4(P)_M \to H_2(\pi;\cy 2)$$
    on an element $[N,f]$ by the projection of the difference $[N,f] - [M,c]$ from 
$$\ker\,(\Omega^{Spin}_4(P) \to H_4(P;\bZ))$$ to the subquotient
$E^{2,2}_\infty =  H_2(\pi;\cy 2)$ in the Atiyah-Hirzebruch  spectral sequence. 
\begin{lemma}
$\rho_M \colon \Omega_4(P)_M \to H_2(\pi;\cy 2)$ is an injection.
\end{lemma}
\begin{proof}
This follows from Proposition \ref{prop:bordism} and the definition of the subset
$\Omega_4(P)_M$.
\end{proof}
Similarly, we define $\Omega_4(M)_M$ and obtain a bijection 
$$\hat\rho_M\colon \Omega_4(M)_M \xrightarrow{\approx} H_2(M;\cy 2),$$ since elements of $\Omega_4(M)_M$ are represented by degree one normal maps. The map 
$$\hat\theta_M\colon \Omega_4(M)_M \to \tilde L_4(\Zpi)$$
is defined by $\hat\theta_M(N, g) = \theta_M(N, c\circ g)$, where  $g\colon N \to M$ represents a bordism element in $\Omega_4(M)_M$. Recall that there is a ``universal" assembly map homomorphism 
$$\kappa_2\colon H_2(\pi; \cy 2) \to  L_4(\Zpi)$$ defined
for any group (see \cite[\S 3]{Kirby:2001}). We have a version of \cite[Lemma 5.11]{HKT09} in our setting.
\begin{lemma} \label{lem:tenfive}
$\theta_M = \kappa_2\circ \rho_M$.
\end{lemma}
\begin{proof}
The following diagram commutes:
$$\xymatrix{&&  L_4(\Zpi)\cr
\Omega_4(M)_M\ar[r]^{c_*}\ar@/^2pc/[urr]^{\hat\theta_M} \ar[d]^{\hat\rho_M}&\Omega_4(P)_M \ar[ur]^{\theta_M}\ar[d]^{\rho_M} && \cr
H_2(M;\cy 2)\ar[r]^{c_*} &H_2(\pi;\cy 2)\ar@/_2pc/[uur]^{\kappa_2}&
}$$
 The outer composition $\hat\theta_M = \kappa_2\circ c_*\circ \hat\rho_M$ holds by the same argument given in the proof of \cite[Lemma 5.11]{HKT09}. The elements in $\Omega_4(M)_M$ are represented by degree 1 normal maps $(N,g)$ covered by a bundle map $\nu_N \to \nu_M$, since  $g^*(\nu_M) \cong \nu_N$ by \cite{Dold:1959,kirby-siebenmann1}.   The required formula now follows from Wall's characteristic class formula for surgery obstructions (see \cite{Davis:2005}).

The map $c_*\colon H_2(M;\cy 2) \to H_2(\pi;\cy 2) $ is surjective, $\hat\rho_M$ is a bijection, and $\rho_M$ is an injection. Hence the formula 
$\theta_M = \kappa_2\circ \rho_M$ follows from the commutivity of the inner square.
\end{proof}
\begin{corollary}\label{cor:fundclass} Let $M$ be a \spinp manifold, and let $P = P(M)$.
Suppose that $[N, f]\in \Omega^{Spin}_4(P)$, with 
$f$ a reduced $3$-equivalence and $\sign(N) = \sign(M)$.
 If $f_*[N] = c_*[M]$ and the map $\kappa_2\colon H_2(\pi;\cy 2) \to L_4(\Zpi)$ is injective, then $[N,f] = [M,c] \in \Omega^{Spin}_4(P)$.
\end{corollary}
\begin{proof}
This is a version of \cite[Corollary 5.12]{HKT09} and the proof is analogous. Since $f$ is a reduced $3$-equivalence and  $f_*[N] = c_*[M]$, Corollary \ref{cor:reducediso} implies that $h_N$ and $h_M$ are isometric (via the given reference maps), and hence $(N,f)$ is  a normal structure. The bordism group $\Omega^{Spin}_4(P)$ is detected by Proposition \ref{prop:bordism}, and the difference $[M,c] -[N,f]$ projects to zero in $H_4(P;\bbZ)$ since $f_*[N] = c_*[M]$. By definition,  the map $\rho_M(N,f)$ is the projection of the difference $[M,c] -[N,f]$  to the subquotient $H_2(\pi;\cy 2)$. Since $f$ is a reduced $3$-equivalence, $\theta_M(N,f) = 0$ by Lemma \ref{lem:tenthree}, and since $\kappa_2$ is injective, $\rho_M(N,f)$ must be zero by Lemma \ref{lem:tenfive}. Furthermore, since $\sign(N) = \sign(M)$, the  elements $[N,f]$ and $[M,c]$ are bordant in $\Omega_4(P)$. 
\end{proof}

\section{Classification of \spinp $4$-manifolds with $\cd \pi _1(M)\leq 3$  }\label{sec:eleven}

In this section we state and prove our main classification result, Theorem \ref{thm:maintame},  for \spinp topological $4$-manifolds with $\cd \pi_1(M) \leq 3$. We conclude by proving Theorem A and Theorem B.

For a finitely presented group $\pi$ with $\cd\pi \leq 3$, let $b_3(\pi)$ denote the minimum number of generators for $H^3(\pi;\La)$ as a $\La$-module. Note that $b_3(\pi)$ is bounded above by the $\La$-rank of $C_3(L)$, where $L$ is a minimal aspherical $3$-complex with $\pi_1(L) = \pi$. We recall a definition from \cite{HKT09}.
\begin{definition}\label{def:WAA}
A group $\pi$ satisfies properties (W-AA) whenever
\begin{enumerate}
\item The Whitehead group $\wh(\pi)$ vanishes. 
\item The assembly map $A_5\colon H_5(\pi;\bbL_0) \to L_5(\Zpi)$ is surjective. 
\item The assembly map $A_4\colon H_4(\pi;\bbL_0) \to L_4(\Zpi)$ is injective. 
\end{enumerate}
\end{definition}
These properties hold whenever the group $\pi$ is torsion-free   and satisfies the Farrell-Jones isomorphism conjectures in $K$-theory and $L$-theory.
These conjectures have been verified for many classes of groups, and in particular for all \RAAGs\ 
 (see \cite{Bartels:2012}, \cite{Bartels:2014}).

\begin{theorem}\label{thm:maintame}
 Let $\pi$ be a finitely presented group with   $\cd \pi \leq 3$ satisfying the properties \textup{(W-AA)}.
 If $M$ and $N$ are closed, oriented, \spinp, TOP $4$-manifolds with fundamental group $\pi$, then any isometry between the quadratic $2$-types of $M$ and $N$ is stably realized by an $s$-cobordism between $M \Sharp r(S^2 \times S^2)$ and $N \Sharp r(S^2 \times S^2)$, for $r \geq b_3(\pi)$.
\end{theorem}
\begin{remark}
Note that 
 the main theorem of \cite{HKT09} applies to groups $\pi$ with $\cd \pi =2$, unless $\pi$ has geometric dimension $3$. Any such group would be a counterexample to  the famous Eilenberg-Ganea conjecture.
\end{remark}
We divide the proof of Theorem \ref{thm:maintame} into the following steps. 
\begin{enumerate}

\item  
Reduced $3$-equivalences $c_M\colon M \to P$ and $c_N\colon N \to P$ arising from an isometry $(\alpha, \beta) \colon Q(M) \cong Q(N)$ of quadratic $2$-types satisfy
$$(c_M)_*[M] = (c_N)_*[N]  \in H_4(P;\bZ).$$
 This is the corresponding result to \cite[Theorem 5.13]{HKT09}, and the proof follows a similar outline (see details below).  Note however that \cite[Lemma 5.16]{HKT09} intended to cite the paper of Whitehead \cite[p.~62]{whitehead2} for the properties used of the $\Gamma$-functor.
\item  The assumption in (W-AA) that the assembly map
$A_4\colon H_4(\pi;\bbL_0) \to L_4(\Zpi)$ is injective implies that the map $\kappa_2\colon H_2(\pi;\cy 2) \to L_4(\Zpi)$ is injective.
\item Suppose that $M$ and $N$ are closed, oriented, \spinp, topological (or smooth)   $4$-manifolds. If $M$ and $N$ have isometric quadratic $2$-types, then there are reduced normal $2$-smoothings $M \to P$ and $N \to P$ which are bordant in $\Omega^{Spin}_4(P)$. This follows as in \cite[Corollary 5.14]{HKT09} adapted to our setting: the details are given in Proposition \ref{prop:seventhree} and  Corollary \ref{cor:fundclass}.
The assumption that $M$ and $N$ have even reduced $w_2$-type is used essentially in this step (see Section \ref{sec:ten}).
\item We show how to apply \cite[Theorem 4, p.~735]{Kreck99} of Kreck's modified surgery theory to obtain an $s$-cobordism between $M$ and $N$,
after a specified stabilization. 
\end{enumerate}

For step (i),  we recall that the key ingredient in the proof of \cite[Theorem 5.13]{HKT09} is the injectivity of the 
  cap product map:
 $$ \omega_P\colon H_4(P;\bZ) \to \Hom_\La (H^2(P;\La), H_2(P;\La)). $$
 defined by cap product with the class $\tr(\alpha) \in H^{lf}_4(\widetilde P; \bZ)^\pi$, 
 for any class $\alpha \in H_4(P;\bZ) $. The formula
 $$b_\alpha(u,v)= \la u  \cup v, \tr(\alpha) \ra, \quad \text{for\ } u,v \in 
 H^2(P;\La),$$
 defines a  symmetric bilinear cup product form
  $$H^2(P;\La) \times H^2(P; \La) \to \bZ,$$
  with $\ad b_\alpha \colon H^2(P;\La) \to H^2(P;\La)^*$ given by cap product with $\tr(\alpha)$, followed by the evaluation map
   $\hat e\colon H_2(P;\La) \to H^2(P;\La)^*$. Since $H_2(P;\La) \cong \pi_2(P)$ is a stably-free $\La$-module, we have the isomorphisms:
  $$\Hom_\La (H^2(P;\La), H_2(P;\La)) \xleftarrow{\approx} \Hom_\La (\pi_2(P)^*, \pi_2(P))$$
 induced by the dual of the evaluation map  $e\colon H^2(P;\La) \to H_2(P;\La)^* $ and the calculations in the proof of Lemma \ref{s_M star}. The two evaluation maps are related by the formula $\hat e = e^* \circ \gamma$, where 
 $\gamma\colon H_2(P;\La) \to H_2(P;\La)^{**}$ is the ``double-dual" isomorphism. We have a natural inclusion
 $$\Gamma(\pi_2(P))^\pi  \subset \Hom_\La (\pi_2(P)^*, \pi_2(P))$$
 by the properties of the $\Gamma$-functor (see \cite[p.~62]{whitehead2}, \cite[p.~144]{HKT09}).  
  
  \medskip
  The same construction applied to the fundamental class $[M] \in H_4(M;\bZ)$ defines the $\pi$-invariant cup product form
  $b_M (u,v) = \varepsilon_1(s_M(Du, Dv))$, for $u,v \in H^2(M;\La)$, where $\varepsilon_1\colon \La \to \bZ$ gives the coefficient at the identity element of $\pi$ and $D\colon H^2(M;\La) \to H_2(M;\La)$ is the Poincar\'e duality isomorphism.
  We note that there exists an isometry $s_M \cong s_N$ of equivariant intersection forms if and only if $b_M$ and $b_N$ are equivariantly isometric.
 
 \medskip
 The map $\omega_P$ with be applied to the image  $( c_{M})_*[M]$ of the fundamental class, whose cap product is given by the composite
 \eqncount
\begin{equation}
\vcenter{\xymatrix@R+5pt@C+15pt{
H^2(M;\La) \ar[r]^{\cap\, \tr [M]}_{\approx} & H_2(M;\La)\ar[d]^{c_*} \\
H^2(P;\La)\ar[u]^{c^*} \ar[r]^{\cap\,  c_*\tr [M]} & H_2(P;\La),}}
\end{equation}
and similarly for $\omega_P$ applied to $(c_{N})_*[N]$.
It follows that there is a commutative diagram (with $\La$-coefficients understood):
 \eqncount
 \begin{equation}\label{eq:omega}
 \vcenter{\xymatrix@C+20pt{
 H_4(M ) \oplus H_4(N ) \ar[d]\ar[r]^(0.3){\omega_M \oplus\, \omega_N}&
 \Hom_\La(H^2(M ), H_2(M )) \oplus \Hom_\La(H^2(N ), H_2(N ))
 \ar[d]\\
 H_4(P;\bZ) \ar[r]^{\omega_P} &\Hom_\La (H^2(P ), H_2(P ))
 }}
 \end{equation}
 where the vertical maps are the sum of the natural maps induced by the reference maps $c^*$ and $c_*$ for $M$ and $N$ respectively. 
 If $M$ and $N$ have isometric quadratic $2$-types, then $b_M \cong b_N$ over $P$ and the difference
 element 
 $$(c_{M})_*[M] - (c_{N})_*[N] \in H_4(P;\La)$$
 maps to zero under $\omega_P$ by the diagram. Therefore  
 $(c_{M})_*[M] = (c_{N})_*[N]$ by the injectivity of $\omega_P$. This completes our outline of step (i).
   
   \medskip
By applying steps (i)-(ii) and Corollary \ref{cor:fundclass}, we may assume that step (iii) is completed. We note that the argument applied to smooth $4$-manifolds produces a smooth spin bordism over $P$. The final step (iv) only works in the topological category.  

The difficulty in step (iv)    is that our reduced normal $2$-smoothings are not given by $2$-connected reference maps $M \to P$ and $N \to P$, so the modified surgery result does not apply directly. 
We now show how a limited amount of stabilization can be used to get  $2$-connected reference maps, and thus complete the proof of Theorem \ref{thm:maintame}.

\medskip
First we construct an abstract diagram using an exact sequence
$$0 \to A \to B \xrightarrow{g} C\to D\to 0$$
and a factorization $B \xrightarrow{j} V \to C$ of the map $g\colon B \to C$.
\eqncount
\begin{equation}\label{eq:diagram}
\vcenter{\xymatrix@R-5pt@C-5pt{
& 0 \ar[d] & 0 \ar[d] & 0 \ar[d] &\\
0 \ar[r] & A \ar[r]\ar[d] & B \ar[r]^{j}\ar[d] \ar[dr]^{g}& V \ar[r]\ar[d] & 0\\
0 \ar[r] & K \ar[r]\ar[d] & B\oplus F \ar[r]^{f}\ar[d] & C \ar[r]\ar[d] & 0\\
0 \ar[r] & E\ar[r]\ar[d] & F \ar[r]^{\bar\phi}\ar[d] \ar[ur]^{\phi}& D \ar[r]\ar[d] & 0\\
& 0  & 0  & 0 &\\
}}
\end{equation}
Let $f(b, x) = g(b) + \phi(x)$, for $b \in B$ and $x\in F$. The map $\phi\colon F \to C$ is a lifting of a surjective map $\bar\phi\colon F \to D$ from a free $\La$-module $F \cong \La^r$. Let $E = \ker\bar\phi$. 
We will apply this diagram to  the universal coefficient sequence
$$0 \to  H^2(\pi; \La) \to  \pi_2(M) \to \pi_2(M)^*\to  H^3(\pi;\La)\to 0$$ from   \eqref{4term for P}, so that $A = H^2(\pi;\La)$, $B = \pi_2(M)$,  $C = \pi_2(P)$ and $D = H^3(\pi;\La)$. The map $g= \ad s_M$, and in our application $\bar\phi\colon F \to H^3(\pi;\La)$ will be given by a (minimal) set of $r $ generators for $H^3(\pi;\La)$ as a $\La$-module. We denote the lifting $\phi$ by $\phi_M\colon F\to \pi_2(P)$ and the resulting map $f$ by $f_M = \ad s_M + \phi_M$.

\begin{remark}    Let $\beta\colon \pi_2(M) \to \pi_2(N)$ be an isometry of the equivariant intersection forms arising from our assumption that $Q(M) \cong Q(N)$. We identify $\pi_1(M)= \pi_1(N) = \pi$. Then 
$$\beta^* \circ \ad s_N \circ \beta = \ad s_M,$$
 and we have reference maps  $(c_M)_* = \ad s_M$, $(c_N)_* = \beta^* \circ \ad s_N$
 as in the proof of Proposition \ref{prop:seventhree}. It follows that $\beta^*= (c_M)^*\circ ((c_N)^*)^{-1}\colon \pi_2(N)^* \to \pi_2(M)^*$ induces an isometry 
 $$\beta^{**} \circ \ad h_M \circ \beta^* = \ad h_N$$ of the reduced intersection forms. After choosing the map $\phi_M \colon F \to \pi_2(P) = \pi_2(M)^*$, we define 
 $\phi_N = (\beta^*)^{-1} \circ \phi_M\colon F \to \pi_2(N)^*$.
  It follows that the pull-back forms $(\phi_M)^*(h_M) = (\phi_N)^*(h_N)$ are equal. The remaining map $f_N \colon \pi_2(N) \oplus F \to \pi_2(N)^*$ used in diagram \eqref{eq:diagram} is defined by 
  $f_N = \ad s_N + \phi_N$.
\end{remark}

\begin{proof}[The proof of Theorem \ref{thm:maintame}]
Let $(W, \Phi)$ be a spin bordism over $P$ from $(M, c_M)$ to $(N, c_N)$, with reference map $\Phi\colon W \to P$. Since the reference maps $c_M$ and $c_N$ on the boundary of $W$ are not $2$-connected (whenever $H^3(\pi;\La) \neq 0$), we 
modify our problem by applying the algebra above to our geometric setting.  We form 
$$M' = M \Sharp r(S^2 \times S^2)$$ 
by performing surgery on null-homotopic circles in $M$.  
Let $\theta_M = -(\phi_M)^*(h_M)$ denote  the hermitian form on $F$ pulled back from 
\emph{minus} the reduced intersection form 
$$h_M \colon \LM \times \LM \to \La$$
 on $\LM = \pi_2(M)^*$,  so that the map $\phi_M$ is an isometry with respect to $-h_M$. Since $w_M = 0$ and $F$ is a free $\La$-module, we may write $\theta_M = \lambda + \lambda^*$ for some sesquilinear form $(F, \lambda)$. Let $\{a_1, \dots, a_r\}$ denote a base for $F$, and $\{b_1, \dots, b_r\}$ denote the dual base for $F^*$. 

Define a hermitian form $(F \oplus F^* , \psi_M)$ 
by the formulas 
$\psi_M(a_i, a_j) = \theta_M(a_i, a_j)$ and $\psi_M(a_i, b_j)  = \delta_{ij}$,   $\psi_M(b_i, b_j) = 0$. 
Let $\{e_1, \dots , e_r,  f_1, \dots, f_r\}$ denote the standard hyperbolic base for $H(\La^r)$, and define an isometry
$$k\colon (F \oplus F^* ,\psi_M) \to H(\La^r)$$
by the explicit formulas $k(a_i) = e_i + \sum_j \alpha_{ij} f_j$, where $\lambda =(\alpha_{ij})$ in matrix form, and $k(b_i) = f_i$. Let $k_0 = p_F \circ k^{-1}$, where $p_F\colon F\oplus F^* \to F$ is the first factor projection.
The reference map $c_M\colon M \to P$,  with induced map
$$ (c_M)_*\colon \pi_2(M) \xrightarrow{\ad s_M} \pi_2(M)^* = \pi_2(P)$$
 can be extended to a $2$-connected reference map $c_{M'}\colon M' \to P$, with induced map
 \eqncount
\begin{equation}
(c_{M'})_*\colon \pi_2(M') = \pi_2(M) \oplus H(\La^r) \xrightarrow{\id \oplus k_0} \pi_2(M) \oplus F \xrightarrow{f_M} \pi_2(P),
\end{equation}
using the map $f_M\colon \pi_2(M) \oplus F \to \pi_2(P)$ from  diagram \eqref{eq:diagram}. Recall that $B = \pi_2(M)$, $C = \pi_2(P)$ and $f_M = \ad s_M + \phi_M $.
The trace of the surgeries provides a spin bordism $(Y, \varphi_M)$ from $(M, c_M)$ to $(M',c_{M'})$   over $P$, where $\varphi_M\colon Y \to P$ is  $2$-connected. 

By the same construction,  the given reference map $c_N\colon N \to P$
  may be extended to a $2$-connected reference map $c_{N'}\colon N' \to P$, 
  where $N' = N \Sharp r(S^2 \times S^2)$, 
 with induced map
\begin{equation*}
(c_{N'})_*\colon \pi_2(N') = \pi_2(N) \oplus H(\La^r) \xrightarrow{\id \oplus k_0} \pi_2(N) \oplus F \xrightarrow{f_N} \pi_2(N)^* \xrightarrow{\beta^*} \pi_2(P),
\end{equation*}
 using the  corresponding $f_N$ from diagram \eqref{eq:diagram}. 
Note that since $\phi_N = (\beta^*)^{-1} \circ \phi_M$, we have $\theta_N = -(\phi_N)^* = \theta_M$ and $\psi_N = \psi_M$.

We obtain a spin bordism $(Z, \varphi_N)$ from  $(N', c_{N'})$  to $(N, c_N)$ over $P$, with reference map $\varphi_N\colon Z \to P$. Now the union of these three bordisms
$$(W', \Phi') = (-Y \cup W \cup Z, \varphi _M\cup \Phi \cup \varphi_N)$$
provides a spin bordism from $(M', c_{M'})$ to  $(N', c_{N'})$, with $2$-connected reference maps $c_{M'}$ and $c_{N'}$ on the two boundary components.

\begin{definition}
Let $P'$ denote the fibration over $K(\pi,1)$ defined by $\pi_2(P') = \pi_2(P) \oplus F^*$, with $k$-invariant pushed forward from $P$ via the inclusion $\pi_2(P ) \to \pi_2(P')$ onto the first summand. 
\end{definition}
Note that we also have a map $P' \to P$ induced by the the first factor projection 
$$p_1\colon \pi_2(P') = \pi_2(P) \oplus F^* \to \pi_2(P).$$
The reference map $c_M\colon M \to P$,  with induced map
$$ (c_M)_*\colon \pi_2(M) \xrightarrow{\ad s_M} \pi_2(M)^* = \pi_2(P)$$
 can be extended to a $2$-connected reference map $\hat c_{M'}\colon M' \to P'$, with induced map
 \eqncount
\begin{equation}
(\hat c_{M'})_*\colon \pi_2(M') = \pi_2(M) \oplus H(\La^r) \xrightarrow{id\, \oplus k^{-1}} \pi_2(M) \oplus F \oplus F^* \xrightarrow{f_M\oplus\,  id_{F^*}} \pi_2(P)\oplus F^*
\end{equation}
with the property that $(c_{M'})_* = p_1\circ (\hat c_{M'})_*$.

\medskip
The same construction will be applied to $(N', c_{N'})$ to produce a $2$-connected reference map to $P'$. The next step is to compare  $(M', \hat c_{M'})$ and $(N', \hat c_{N'})$.

\medskip
\noindent
{\bf Claim 1}: \emph{The reference maps $(M', \hat c_{M'})$ and $(N', \hat c_{N'})$ induce an isometry $s_{M'} \cong s_{N'}$ of equivariant intersection forms over $P'$.}

\begin{proof} We check directly from the definitions that the isometry
$$\beta\oplus \id \colon (\pi_2(M), s_{M})\oplus H(\La^r) \xrightarrow {\approx} (\pi_2(N), s_{N})\oplus H(\La^r)$$
arising from a given isometry $(\id, \beta)\colon Q(M) \cong Q(N)$ of quadratic $2$-types is compatible with the reference maps to $P'$, as defined above. Since $\psi_M = \psi_N$, this reduces to checking the equation 
$$f_M = \beta^*  \circ f_N \circ (\beta \oplus \id) \colon \pi_2(M) \oplus F \to \pi_2(P),$$
by verifying that
\begin{align*}
(\beta^* \circ f_N \circ (\beta\oplus \id))(b, x) &= (\beta^*\circ f_N)(\beta(b), x) = \beta^*( \ad s_N (\beta(b)) +\beta^*(\phi_N(x))\\ &= \ad s_M (b) + \phi_M(x) = f_M(b,x)
\end{align*}
for all $b \in \pi_2(M)$ and all $x \in F$.
\end{proof}

\medskip
\noindent
{\bf Claim 2}: \emph{The elements $(M', \hat c_{M'})$ and $(N', \hat c_{N'})$ are spin bordant over $P'$.}
 
 \begin{proof}
We already have 
 a spin bordism $(W',\Phi')$ from $(M', c_{M'})$ to $(N', c_{N'})$ over $P$, so the difference element 
 $$\alpha:=  [(N', \hat c_{N'})] - [(M', \hat c_{M'})]$$
 is in the kernel of the natural map  $\Omega_4^{Spin}(P') \to \Omega_4^{Spin}(P)$ induced by the map $P' \to P$ defined above. By the Atiyah-Hirzebruch spectral sequence (as in Section \ref{sec:ten}), we obtain an injection
 $$\Omega_4^{Spin}(P') \to \Omega_4^{Spin}(P) \oplus H_4(P';\bZ)$$
 and the image of $\alpha \in \Omega_4^{Spin}(P')$ in $H_4(P';\bZ)$ is determined by the difference of the fundamental classes
  $(\hat c_{N'})_*[N'] - (\hat c_{M'})_*[M'] \in H_4(P';\bZ)$. 
  
  \medskip
  By  Claim 1, we have an isometry $s_{M'} \cong s_{N'}$ over $P'$, and hence an equivariant isometry $b_{M'} \cong  b_{N'}$ over $P'$. Moreover, we have the analogue of diagram \eqref{eq:omega} for $M'$, $N'$ and $P'$, and the map $\omega_{P'}$ is again injective. The same argument used in our outline of step (i) now applies to show that 
   $(\hat c_{N'})_*[N'] = (\hat c_{M'})_*[M'] \in  H_4(P';\bZ)$, as required. 
 \end{proof}

We now have a spin bordism $(\widehat W', \widehat \Phi')$ from $(M', \hat c_{M'})$ to  $(N', \hat c_{N'})$, with $2$-connected reference maps $\hat c_{M'}$ and $\hat c_{N'}$ on the two boundary components. After further surgeries relative to the boundary, we may assume that $\widehat\Phi'\colon \widehat W' \to P'$ is also 2-connected.
Diagram \eqref{eq:diagram} also determines the structure of
$$KH_2(M') := \ker((\hat c_{M'})_*: \pi_2(M') \to \pi_2(P')).$$
Since the free module $F \cong (\La^r\times \{0\})$ is required to map surjectively onto $H^3(\pi;\La)$, we may take $r = b_3(\pi)$.

\medskip
\noindent
{\bf Claim 3}: \emph{The form $s_{M'}$ restricted to  $\ker (\hat c_{M'})_*\subset \pi_2(M')$ is identically zero.}

\begin{proof}
Since $\ker (\hat c_{M'})_* \cong K = \ker f_M$ (from Diagram \ref{eq:diagram}), 
it is enough to show that the form on $K = \ker f_M$ induced  by 
$$(\pi_2(M'), s_{M'}) \cong (\pi_2(M), s_M)  \oplus (F\oplus F^*, \psi_M)$$
 is identically zero.
To check this, we compute the hermitian form $s_M \oplus \psi_M \cong s_{M'}$ on elements $z = (b, x,0)$ and $z'=(b', x',0)$ of $K\subset \pi_2(M) \oplus F\oplus F^*$ via
\begin{align*}
(s_M \oplus \psi_M)(z, z') = s_M(b, b') + \theta_M(x, x') &=   s_M(b, b') - h_M(\phi_M(x), \phi_M(x'))\\
 &=  s_M(b, b') - \ad h_M (\phi_M(x))(\phi_M(x')).
\end{align*}
But $\ad s_M(b) = -\phi_M(x)$ and $\ad s_M(b') = -\phi_M(x')$,  since our elements lie in $K = \ker f$, so we obtain
$$(s_M \oplus \psi_M)(z, z')=   s_M(b, b') - \ad h_M (\ad s_M(b))(\ad s_M(b')) = s_M(b, b') -
\la \gamma(b), \ad s_M(b')\ra = 0$$
by the formula in \eqref{eq:gamma1}.
\end{proof}

We now have the setting to apply Kreck \cite[Theorem 4, p.~735]{Kreck99}, which measures the obstruction to surgery on $(\widehat W', \widehat\Phi')$, relative to the boundary, to obtain an $s$-cobordism between $(M', \hat c_{M'})$ and $(N', \hat c_{N'})$.    Claim 3 above  implies that the modified surgery obstruction is zero, by \cite[Prop.~8, p.~739]{Kreck99}. The role of the first two properties in (W-AA) is explained in \cite[Theorem 2.6]{HKT09}, showing that we may modify our bordism $(\widehat W', \widehat \Phi')$ by the action of $L_5(\Zpi)$  to obtain an $s$-cobordism. This completes step (iv) and the proof of Theorem \ref{thm:maintame}. 
\end{proof}

\begin{proof}[The proof of Theorem A]
For $\pi$ a \RAAG\ with $\cd\pi \leq 3$, the conditions (W-AA) were established by
 Bartels and L\"uck \cite{Bartels:2012}. This shows that Theorem A follows from
 Theorem \ref{thm:maintame}.
\end{proof}

\begin{proof}[The proof of Theorem B]
Suppose that $M$ and $N$ are topological manifolds as in the statement, with $s_M$ stably isometric to $s_N$. Then by Proposition \ref{prop:movek}, we conclude that $Q(M)$ is stably isometric to $Q(N)$. Hence, after stabilization by sufficiently many copies of $S^2 \times S^2$, we may apply all but the final step in  the proof of  Theorem \ref{thm:maintame} to conclude that $M$ and $N$  are spin bordant over $P$,  and hence stably homeomorphic. For these steps in the argument we only need the injectivity of $\kappa_2$,  which is property (W-AA)(iii).  In the smooth case, we apply essentially the same arguments above to show that $M$ and $N$ are smoothly spin bordant over the reduced type $P \times BSPIN$, and it follows that $M$ and $N$ are stably diffeomorphic by \cite[Theorem C]{Kreck99}.
\end{proof}
\section{Appendix: the hypercohomology spectral sequence}\label{sec:append}

We include some background material and references for homological algebra for the reader's convenience. 

\nr{1} The universal coefficient spectral sequence is a special case of of the \textbf{hypercohomology spectral sequence} defined in Benson \cite[Proposition 3.4.3]{benson2}, Another good reference is Benson and Carlson \cite[\S2 and \S4]{Benson:1994}. If $\bC$ and $\bD$ are chain complexes over a ring $\La$, with $\bC$ bounded above and $\bD$ bounded below, then there is a spectral sequence
$$ \Ext^p_{\La}(H_q(\bC), \bD) \Rightarrow \Ext^{p+q}_{\La}(\bC, \bD)$$
converging to the hypercohomology groups (see Benson and Carlson \cite[\S 2.7]{benson1} for the definition of $\Ext$ for chain complexes). This is a first quadrant spectral sequence, so the indexing convention above is standard ($p$ runs along the $x$-axis). The differential $d_r$ has bi-degree $(r, -r+1)$. 

\nrb{1a}  For $X$ a connected finite CW-complex with $\pi_1(X) = \pi$ and $\La = \Zpi$, we can take $\bC = C(X;\La)$ and $\bD = \La$ (a $0$-dim chain complex with $\La$ in degree zero). Then we have the \textbf{universal coefficient spectral sequence}
$$ \Ext^p_{\La}(H_q(X;\La)), \La) \Rightarrow \Ext^{p+q}_{\La}(C(X;\La), \La) = H^{p+q}(X;\La).$$
For $q=0$, the terms along the $x$-axis are just $\Ext^p_{\La}(\bZ, \La) = H^p(\pi;\La)$. On the $y$-axis we have $\Hom_{\La}(H_q(X; \La), \La)$. To compute $H^2(X;\La)$ for $X$ a closed, oriented $4$-manifold, the spectral sequence produces the 4-term exact sequence:
$$ 0 \to H^2(\pi;\La) \to H^2(X;\La) \to \Hom_{\La}(H_2(X; \La), \La) \to H^3(\pi;\La) \to 0$$ 
that we use throughout the paper. Note that the map to $H^3(\pi;\La)$ is a $d_2$ differential, which is surjective since $H^3(X;\La) = H_1(X;\La) = H_1(\wX; \bZ) = 0$.

\nrb{1b} There are two usual definitions of the $\Ext$-functors, namely in terms of projective resolutions (the usual homological algebra), or in terms of \emph{multiple extensions} due to
MacLane \cite[Chap.~III, \S\S 5-6]{MacLane:1967}. The multiple extensions version is particularly useful in describing connecting maps in long exact sequences (see Benson and Carlson \cite[\S 2.6]{benson1}). 

\nr{2} \textbf{A formula for the $d_2$ differential}:  Let $X$ be a finite CW complex with $\pi:=\pi_1(X)$ and integral group ring $\Lambda = {\bZ}\pi$. Let $A$ be a right $\La$-module.
We use the universal coefficient spectral sequence which has $E_2$-page 
$$E^{p,q}_2 = \Ext^p_\La(H_q(X;\La),A).$$ 
to calculate the cohomology $H^*(X; A)$.
The $d_2$ differential 

$$d_2\colon E^{p,q}_2 = \Ext^p_\La(H_q(X;\La),A) \to 
E^{p+2,q-1}_2 = \Ext^{p+2}_\La(H_{q-1}(X;\La), A) $$
is given in terms of  the multiple extensions description for $\Ext$ by the following ``splicing" formula. 

\nrb{2a}
Let $C_* : = C(X;\La)$ denote the chain complex of $X$ with $\La$-coefficients, and consider the usual sequences
of cycles $Z_*:= Z_*(C)$ and boundaries $B_*:=B_*(C)$
\begin{align*}
& 0 \to B_{q+1} \to Z_{q}  \to H_q \to 0\\
& 0 \to Z_{q} \to C_{q}  \to B_q \to 0\\
& 0 \to B_{q} \to Z_{q-1}  \to H_{q-1}\to 0
\end{align*}
 for computing $H_* := H_*(X;\La)$. We can put these together to obtain the
 $4$-term  exact sequence
\eqncount
\begin{equation}\label{one}
 0 \to H_q \to C_q/B_{q+1} \to Z_{q-1} \to H_{q-1} \to 0.
 \end{equation}

\noindent
To define $d_2$, we represent an element $\alpha \in  \Ext^p_\La(H_q,A)$
by a multiple extension
\eqncount
\begin{equation}\label{two}
 0 \to A \to S_{p-1} \to S_{p-2} \to \dots \to S_0 \to H_q\to 0
 \end{equation}
as in the description of $\Ext$ given by MacLane \cite[Chap.~III, Theorem 6.4]{MacLane:1967}, and splice the resolutions \eqref{one} and \eqref{two} together to get a representative for 
$$d_2(\alpha) \in \Ext^{p+2}_\La(H_{q-1};A). $$

\nr{3} Finally we relate this description of  the $d_2$ differential as a Yoneda splice to the usual definition.
In general (even when $C_*$ is not a projective chain complex),  the hypercohomology spectral sequence is defined (see Benson \cite[\S 3.4]{benson2}) as a spectral sequence
of a double complex $$E_0 ^{p,q}= \Hom _{\Lambda } ( P_{p,q}, A)$$ where for each $q\geq 0$, $P_{*, q} \to C_q$ is a projective resolution
of $C_q$, and the chain maps $P_{*, q} \to P_{*, q-1}$ are induced by the boundary maps $C_q \to C_{q-1}$. Applying the horseshoe lemma
to the extensions $$ 0 \to B_{q+1} \to Z_q \to H_q \to 0$$ and   $$0 \to Z_q \to C_q \to B_q \to 0,$$
we can construct a projective resolution for $C_q$ of the form $P_*(B_{q+1}) \oplus P_* (H_q) \oplus P_* (B_q)$, where $P_*(B_{q+1})$,  $P_* (H_q)$, and  $ P_* (B_q)$ are projective resolutions of $B_{q+1}$, $H_q$, and $B_q$, respectively.  The boundary map on this projective resolution is of the form $$ \partial _1(a, b, c) =(\partial a+f(b)+g_1(c), \partial b + g_2(c), \partial c)$$
where $f\colon  P_* (H_q) \to P_* (B_{q+1})$ and $(g_1, g_2)\colon  P_* (B_q) \to P_*(Z_q)=P_* (B_{q+1})\oplus P_* (H_q)$ are degree $-1$ chain maps. Note that $f\circ g_2 = 0$ to check the relation $\partial_1\circ \partial_1 = 0$.
After applying Hom, we will get the horizontal differential $d_1\colon E^{p,q}_0 \to E^{p+1,q}_0$, which can be expressed as the degree $1$ map
$$d_1(a,b,c)=(\delta a, f'(a)+\delta (b), g'_1(a)+g'_2(b)+\delta(c)),$$
 where $f'$ and $g'_1$, and $g'_2$ are duals of the corresponding maps. 

The differential $\partial_0\colon P_* (C_q) \to P_* (C_{q-1}) $ is induced by the maps $C_q\to B_q$ 
and the composite $B_q \to Z_{q-1} \to C_{q-1}$. 
Hence by the construction of the projective resolution $$P_*(C_q)=P_*(B_{q+1}) \oplus P_* (H_q) \oplus P_* (B_q),$$ the induced map is
$\partial_0 (a,b,c)=(c,0,0)$ and its dual gives the vertical differential $d_0\colon E^{*,q-1}_0 \to E^{*,q}_0$. Then one defines
$E_1^{p,q} = H(E_0^{p,q}, d_0)$ and $E_2^{p,q} = H^p(H^q(E_0^{*,*},d_0),d_1)$.

In our case, we have $E_2^{p,q} = \Ext ^p _{\Lambda } (H_q, A)$.  The next differential in the spectral sequence
 $$d_2\colon  \Ext ^p _{\Lambda } (H_q, A) \to \Ext ^{p+2} _{\Lambda } (H_{q-1} , A), $$
 is now defined using the diagram
 $$\xymatrix{&\Hom_\La(P_{p+1}(C_{q-1}),A) \ar[r]^{d_1}\ar[d]^{d_0}&\Hom_\La(P_{p+2}(C_{q-1}),A)\cr
\Hom_\La(P_{p}(C_{q}),A)\ar[r]^{d_1}&\Hom_\La(P_{p+1}(C_{q}),A)&}$$
 as the map obtained by first applying the horizontal $d_1$ differential, then lifting via the vertical $d_0$ differential, and then applying the horizontal $d_1$ again to the lifted element (see Benson \cite[p.~196]{benson2} for the notation).
 
Let $u \in \Hom _{\Lambda} (P_p (H_q) , A )$ be a degree $p$ cocycle representing an $\Ext$-class in $\Ext^p _{\Lambda} (H_q, A)$. 
On $E_0 ^{p,q}$ this is represented by the class $(0, u, 0) \in E_0^{p,q} = \Hom _{\Lambda} (P_p (C_q) , A)$. Applying the horizontal differential on it, we obtain $(0, 0, g'_2 (u))$ in $E_0^{{p+1},q}$, where $g'_2$ is the dual of the chain map $g_2$. 

Since the vertical differential $d_0$ is induced by dualizing $\partial _0$, we see that the element
 $(g'_2(u) , 0, 0) \in E_0 ^{p+1, q-1}$ maps to $(0, 0, g'_2(u))$ under the vertical differential. 
Now by applying the horizontal differential again we get  $(0, f' (g'_2 (u)), g'_1 (g'_2 (u))) \in E^{p+2,q-1}_0$. The element $f'(g'_2(u))=(g_2 f)'(u)$ gives a degree $(p+2)$-cocycle in $\Hom _{\Lambda} ( P_*(H_{q-1} ) , A)$. The cohomology class $[(g_2f)'(u))]$ is the image of $[u]$ under $d_2$.    

Note that $f$ is the chain map
$f\colon  P_* (H_{q-1}) \to P_* (B_{q})$, and $g_2$ is the second coordinate map of the chain map  $(g_1, g_2)\colon P_* (B_{q}) \to P_*(Z_{q})=P_* (B_{q+1})\oplus P_* (H_{q})$. The chain map $f$ is induced by the extension 
$$0 \to B_{q} \to Z_{q-1} \to H_{q-1} \to 0.$$
 In addition, we can regard $g_2$ as the chain map induced by 
 $$0 \to H_{q} \to C_{q}/ B_{q+1} \to B_{q} \to 0.$$  
 The multiple extensions description of the $\Ext$-group shows that 
applying these chain maps corresponds to Yoneda splice with the corresponding extension classes. This completes the identification of the $d_2$ differential. \qed


\providecommand{\bysame}{\leavevmode\hbox to3em{\hrulefill}\thinspace}
\providecommand{\MR}{\relax\ifhmode\unskip\space\fi MR }
\providecommand{\MRhref}[2]{%
  \href{http://www.ams.org/mathscinet-getitem?mr=#1}{#2}
}
\providecommand{\href}[2]{#2}

\end{document}